\documentclass[twoside,11pt,reqno]{amsart}
\usepackage{graphicx,lscape}
\usepackage{amsmath}
\usepackage{amsfonts}
\usepackage{amssymb,mathrsfs}
\usepackage{url,ifthen,color,xcolor}
\usepackage{amsthm}
\newtheorem{thm}{Theorem}[section]

\newtheorem{cor}[thm]{Corollary}
\newtheorem{lem}[thm]{Lemma}

\theoremstyle{definition}
\newtheorem{defn}[thm]{Definition}
\newtheorem{exmp}[thm]{Example}
\newtheorem{alg}[thm]{Algorithm}

\theoremstyle{remark}
\newcommand{\cS}{\mathcal{S}}
\newcommand{\cN}{\mathcal{N}}

\newcommand{\bN}{\mathbb{bN}}
\newcommand{\bR}{\mathbb{R}}
\newcommand{\bQ}{\mathbb{Q}}
\newcommand{\bC}{\mathbb{C}}
\newcommand{\bZ}{\mathbb{Z}}

\newcommand{\sA}{\mathscr{A}}
\newcommand{\sC}{\mathscr{C}}
\newcommand{\sG}{\mathscr{G}}

\newcommand{\sM}{\mathscr{M}}

\newcommand{\CC}{\mathbb{C}}
\newcommand{\CTwo}{\mathbb{C}^2}
\newcommand{\Pt}{\mathbb{P}^2}

\newcommand{\ep}{\varepsilon}

\newcommand{\Henon}{H\'{e}non }
\newcommand{\cD}{\mathcal{D}}
\newcommand{\cB}{\mathcal{B}}
\newcommand{\cR}{\mathcal{R}}
\newcommand{\chrec}{\cR}
\newcommand{\chreceps}{\cR_{\ep}}
\newcommand{\chrecepsp}{\cR_{\ep'}}
\newcommand{\twovec}[2]{\begin{pmatrix} #1 \\ #2  \end{pmatrix}}

\renewcommand{\Re}{\text{Re}}
\renewcommand{\Im}{\text{Im}}

\newcommand{\abs}[1]{\left\arrowvert #1\right\arrowvert}
\newcommand{\norm}[1]{\left\| #1\right\|}
\newcommand{\onorm}[1]{\abs{#1}}
\newcommand{\snorm}[1]{\norm{#1}}
\newcommand{\enorm}[1]{\norm{#1}_e}
\newcommand{\sonorm}[1]{\onorm{#1}}
\newcommand{\inorm}[1]{\norm{#1}_{\text{inf}}}
\newcommand{\supnorm}[1]{\norm{#1}_{\text{sup}}}

\newcommand{\supp}{\text{supp}}

\title[Computability for polynomial skew products]{Computability for Axiom A Polynomial Skew Products of $\CTwo$}

\author[S.~Boyd]{Suzanne Boyd}
\address{Department of Mathematical Sciences\\
University of Wisconsin Milwaukee\\
PO Box 413\\
Milwaukee, WI 53201, 
USA}
\email{sboyd@uwm.edu, ORCID: 0000-0002-9480-4848.}

\author[C. Wolf]{Christian Wolf}
\address{Department of Mathematics and Statistics\\
Mississippi State University\\
Starkville, MS 39759, USA}
\email{cwolf@math.msstate.edu, ORCID: 0000-0002-7976-3574.}
\thanks{C.W. was partially supported by grants from the Simons Foundation (SF-MPS-CGM-637594 and SFI-MPS-TSM-00013897) and PSC-CUNY (TRADB-55-67482, TRADB-56-68594).}

\date{\today}

\begin{document}

\begin{abstract}
    The computability of Julia sets of rational maps on the Riemann sphere has been intensively studied in recent years, see, e.g. \cite{Braverman-Yampolsky-2008, Burr-Wolf-2024} for an overview. For example, by Braverman's results  \cite{Braverman2005, Braverman2006}, hyperbolic and parabolic Julia sets are computable in polynomial time. In this paper, we present the first work on computability related to maps of more than one complex dimension.
    We examine a family of polynomial endomorphisms of $\CTwo$, the polynomial skew products; i.e., maps of the form
    $f(z,w) = (p(z), q(z,w)),$
    where $p$ and $q$ are complex polynomials of the same degree $d\geq 2$.  
    We show that if a polynomial skew product is Axiom A, then its chain recurrent set, which is equal to its non-wandering set and also equal to the closure of the periodic orbits, is computable. Our algorithm also identifies the various hyperbolic sets of different types, i.e., expanding, attracting, and hyperbolic sets of saddle-type.
    One consequence of our results is that Axiom A is a semi-decidable property on the closure of the  Axiom A polynomial skew product locus. Finally, we introduce an algorithm that establishes the lower semi-computability of the hyperbolicity locus of polynomial skew products of a fixed degree. 

\end{abstract}

\maketitle


%%%%%%%%%%%%%%%%%%%%%%%%%%%%%%%%%%%%%%%%%%%%%%%%%%
% for arxiv: 36 pages, 5 figures

%====================================

\section{Introduction}

\subsection{Motivation}

The goal of this paper is to extend computability results for one-dimensional complex dynamical systems to higher-dimensional complex dynamics. Specifically, we apply constructive dynamical methods to develop algorithms which allow us to compute the non-wandering set of Axiom A polynomial skew products of $\bC^2$
to any prescribed level of accuracy.

Computability in dynamical systems has been a topic of significant interest over the past two decades. In particular, the computability/non-computability of Julia sets in one-dimensional complex dynamics (e.g., \cite{BBRY-2011, BBY1, Braverman2005, Braverman-Yampolsky-2008,BY2,D1, DY2018, R2005}) and more recently, the computability of various dynamical invariants including entropy, general invariant sets, spectra, topological pressure, zero-temperature limits, and equilibrium states (\cite{Binder-et-al-2024, BurrWolf2018, BDWY2022, BSW2020, Burr-Wolf-2024, GHR2011, GraccaEtAl2018, HP2025,  HM, HS,HR,Spandl2008}) has attracted significant attention.

One of the most striking findings in the computability theory of one-dimensional complex dynamics is the existence of polynomials with computable coefficients whose Julia sets are non-computable \cite{Braverman-Yampolsky-2008}. Even more surprisingly, the measure of maximal entropy (the Brolin-Lyubich measure which is supported on the Julia set) is always computable \cite{BBRY-2011}. This substantiates a significant distinction between the computability of sets and the associated invariant measures.

A fundamental challenge in computability in dynamical systems stems from the nature of computer algorithms: they require finite input, whereas the mathematical descriptions of various dynamical objects, including invariant sets, may involve infinite data. Consequently, it is not a priori clear whether accurate approximations of such objects are possible. Computability theory addresses this challenge by providing tools to decide whether a finite amount of information is sufficient to compute meaningful approximations of dynamical objects.

In this paper, we answer this question affirmatively for a class of higher-dimensional complex dynamical systems by providing an explicit algorithm that computes the non-wandering set of Axiom A polynomial skew products of $\bC^2$. We obtain rigorous computability results by combining, in a novel way, foundational tools from hyperbolic dynamics, including a quantitative version of the shadowing lemma, graph theory, and various other analytical tools.
We hope that the methods introduced in this paper serve as a foundation for future work on the computability of natural invariant sets for other classes of real and complex dynamical systems.

Because we initiate herein the study of computability for dynamical systems in several complex variables, we provide background material from both computability theory and higher-dimensional complex dynamics, to make the paper more self-contained.

%================================
\subsection{Statement of the Results} 
\label{sec:intro-results}
Let $f:\bC^2\to\bC^2$ be a polynomial skew product, that is, $f(z,w)=(p(z),q(z,w))$, where $p$ and $q$ are complex polynomials with (equal) degree $d\geq 2$. The map $f$ extends to a polynomial endomorphism on two-dimensional complex projective space $P_2(\bC)$ (\cite{Jonsson1999}). Our study exploits the fact that skew products are a natural generalization of polynomials in one dimension, namely $f$ maps the vertical line $\{z\} \times \CC$ to the vertical line $\{p(z)\} \times \CC$, and restricted to a vertical line, $f$ is the polynomial map $w \mapsto q_z(w)=q(z,w)$. We say that $f$ is \textit{Axiom A} if its non-wandering set $\Omega$ is a hyperbolic set and the periodic points are dense in $\Omega$. In this case, the non-wandering set splits into one to four types of compact invariant hyperbolic sets:
\[
\Omega=J_2\cup X_0\cup X_1\cup \Lambda.
\]
The set $J_2\neq \emptyset$ is a uniformly expanding repeller, coinciding with the closure of the repelling periodic points of $f$, and is called the Julia set of $f$. The set $X_0$ is a union of uniformly attracting invariant subsets, which may be empty. The sets $X_1$ and $\Lambda$, which also can be empty, are unions of hyperbolic sets of saddle type, with $\Lambda$'s subsets having expansion in the base (the $z$-coordinate) and contraction in the fibers (the $w$-coordinate), while $X_1$'s sets have the reversed expansion/contraction. 
%. 
We note that for Axiom A polynomial skew products, the non-wandering set $\Omega$ coincides with the chain recurrent set $\cR$ (\cite{Jonsson1999}). Further, $\Omega$ is the union of the chain components of $f$, which is a crucial property in our approach.
We refer to Section~\ref{sec:skewintro} and \cite{Boyd-DeMarco-2008,Jonsson1999} for more discussions of polynomial skew products.

Braverman established  (\cite{Braverman2005,Braverman2006}) that hyperbolic and parabolic Julia sets are computable in polynomial time (also called poly-time computable). Roughly speaking, a subset $C$ of $\bR^\ell$, or $\bC^\ell$, is computable if there exists a Turing machine (a computer program for our purposes) which, on input $n$, outputs a finite set of dyadic points which is $2^{-n}$ close to $C$ in the Hausdorff metric. Intuitively, for planar sets, this means that one can print on a computer screen an accurate picture of the set $C$ at any prescribed error margin. Here, the output of the computer program, finitely many dyadic points, are the coordinates of the pixels visualizing the set $C$. 
Braverman's computability results for Julia sets have been extended to many types of one-dimensional Julia sets. In fact, all but Siegel Julia sets are in a certain sense computable, see \cite{BY2009,Burr-Wolf-2024}. We stress that these computability results generally rely on geometric characterizations of the Julia sets, which are not available in higher dimensions. In particular, all currently existing proofs for poly-time computability depend in one way or another on conformality, which is lacking in two-dimensional complex dynamics.

In this paper, we study the computability, for Axiom A polynomial skew products, of the Julia set $J_2$ as well as the invariant sets which are saddle-type ($X_1$ and $\Lambda$) or attracting ($X_0$). To accomplish this, we develop a novel approach to establish computability based on various tools from hyperbolic dynamics, graph theory, and geometry. We also use the simple but crucial idea to work with a higher iterate of $f$, for which a Turing machine can detect uniform expansion/contraction (as appropriate) for a large enough iterate. Another important feature of our approach is the adaptation of techniques from the work of the first author (under her former name, \cite{Hruska-Skew-2006}) to the realm of computability theory. This adaptation is used to build a model of the dynamics of the map in a neighborhood of its chain recurrent set. The main result of this paper is the following:

\begin{thm}
\label{thm:main}
     Let $f$ be an Axiom A polynomial skew product of $\bC^2$. Then the Julia set $J_2$ of $f$ is computable. Moreover, there is an algorithm that computes all non-empty hyperbolic sets $X_0,X_1$, and $\Lambda$, and determines their emptiness otherwise.
\end{thm}

A specific feature of our algorithm is that, based on the input of the coefficients of the polynomials $p$ and $q$, the program halts if the corresponding map $f$ is Axiom A. As a consequence, we obtain the following:

\begin{cor}
\label{cor:AxiomAsemi-decidable}
    A polynomial skew product of $\bC^2$ being Axiom A is a semi-decidable property on the closure of the Axiom A locus of polynomial skew products of fixed degree. More precisely, there exists a Turing machine that, on oracle access to the coefficients of $p$ and $q$, halts if $f=f_{p,q}$ is Axiom A, and runs forever if $f$ is not Axiom A.
\end{cor}

Another consequence of Theorem \ref{thm:main} is the computability of the measure of maximal entropy of an Axiom A polynomial skew product. Indeed, Binder et.\ al.\ established in \cite{Binder-et-al-2024}  the computability of the measure of maximal entropy (and more generally the computability of equilibrium states of computable H\"older continuous potentials) for distance expanding maps. Applying this result to the Julia set of an Axiom A skew product $f$ we conclude that the unique measure of maximal entropy $\mu_{\rm max}$ of $f$ (which is supported on $J_2$), as well as equilibrium states of computable H\"older continuous potentials on $J_2$, are computable. 

Finally, we apply the hyperbolicity detection feature of our algorithm to derive computability results for the Axiom A locus of polynomial skew products. Given $d\geq 2$ let $\mathscr{A}_d$ denote the space of Axiom A polynomial skew products, for which both $p$ and $q$ are of fixed degree $d$. The set $\mathscr{A}_d$ is an open subset of $\bC^\ell$ for some $\ell\in\bN$ (see \cite{Jonsson1999}, Cor. 8.15 for openness).
An open set $U\subset \bC^\ell$ is lower semi-computable
(also called lower computable) if there exists a Turing machine producing a (possibly infinite) string of points $(x_i,r_i)$, where $x_i$  is a dyadic point in $\bC^\ell$ and $r_i$ is a positive dyadic radius, such that $U = \cup_{i=1}^\infty B(x_i,r_i)$. We obtain the following:

\begin{thm}
\label{thm:LowerComputabilityLocus}Let $d\geq 2$. The locus $\mathscr{A}_d$ of Axiom A polynomial skew products of degree $d$ is lower semi-computable.
\end{thm}

We note that the definition of lower semi-computability does not require the set to be bounded. Indeed, Astorg and Bianchi \cite{Astorg-Bianchi} provide polynomial skew products of fixed degree with unbounded hyperbolic components.

We close this section with the organization of the remaining sections. In Section~\ref{sec:prelim} we provide more detailed background, definitions, theorems, etc., on: computability (subsection~\ref{sec:compute:basic}),
chain recurrence and hyperbolicity (subsection~\ref{sec:dpspaces}), and polynomial skew products (subsection~\ref{sec:skewintro}). In Section~\ref{sec:imagebounds}, we prove a collection of lemmas to enable us to study polynomial skew products from a computability perspective.    
In subsection~\ref{sec:hypatia}, we present our version of the algorithm from \cite{Hruska-Skew-2006}, adapted for computability purposes. 
In Section~\ref{sec:proof-thm-1}, we provide the algorithm for computing the sets $J_2,X_0,X_1,\Lambda$ (any that exist) for an Axiom A polynomial skew product and prove Theorem~\ref{thm:main}. In Section~\ref{sec:results2}, we prove Corollary~\ref{cor:AxiomAsemi-decidable} and Theorem~\ref{thm:LowerComputabilityLocus}.

%====================================
\section{Background/Preliminaries}
\label{sec:prelim}
%====================================

%-----------------------------------
\subsection{Computability}
\label{sec:compute:basic}
%-------------------------------
We are interested in the correctness of computations of certain dynamically defined subsets of $\bR^\ell$, respectively $\bC^\ell$.
Computability theory allows us to guarantee the correctness and accuracy of the computations of these sets.  We recall that a computer can approximate only finitely many real numbers.  Thus, without an accuracy guarantee, a computationally derived approximation of a set could miss interesting features.

We refer to \cite{BBRY-2011, Binder-et-al-2024, Braverman2005, BY2009, BSW2020, Burr-Wolf-2024, GHR2011} 
for more detailed discussions of computability theory.  We use closely related definitions to those in \cite{BY2009} and \cite{BSW2020}. Throughout this discussion, we use a bit-based computation model (this just means information is stored as binary digits), such as a Turing machine (a computer program for our purposes).  One can think of the set of Turing machines as a particular, countable set of functions; we denote $T(x)$ as the output of the Turing machine $T$ based on input $x$.

We start with the definition of computable points in $\ell$-dimensional Euclidean space. 

\begin{defn}
Let $\ell\in \bN$ and $x\in \bR^\ell$. An \emph{oracle} of $x$ is a function $\phi:\bN\to \bQ^\ell$ such that $\Vert \phi(n)-x\Vert < 2^{-n}$. Moreover, we say $x$ is computable
if there is a Turing Machine $T=T(n)$ 
which is an oracle of $x$.
\end{defn} 
It is straightforward to see that rational numbers, algebraic numbers, and some transcendental numbers such as ${ e}$ and $\pi$ are computable real numbers. However, since the collection of Turing machines is countable, most points in $\bR^\ell$ are not computable.
Identifying $\bC^\ell$ with $\bR^{2\ell}$, the notion of computable points naturally extends to $\bC^\ell$.

Next, we define computable functions on Euclidean spaces.
\begin{defn}\label{defcompfunc} Let $D\subset \bR^\ell$. A function $f:D\rightarrow\bR^k$ is {\em computable} if there is a Turing machine $T$ so that for any $x\in D$, any oracle $\phi$ for $x$ and any $n\in \bN, \,\, T(\phi,n)$
is a point in $\bQ^k$ so that $\Vert T(\phi,n)-f(x)\Vert<2^{-n}$.
\end{defn}
We observe that one of the inputs of the Turing machine $T$ in Definition \ref{defcompfunc} is an oracle. Specifically, while the Turing machine $T$ in principle has access to an infinite amount of data, it must be able to decide when the approximation $\phi(m)$ of $x$  is sufficiently accurate to perform the computation of $f(x)$ to precision $2^{-n}$.
We further note that in Definition \ref{defcompfunc}, the input points $x$ are not required to be computable. In fact, any set of $D\subset \bR^\ell$ can be the domain of a computable function $f$. Next, we extend the notion of computable points to more general spaces, called computable metric spaces.

\begin{defn}\label{def:computable}
Let $(X,d_X)$ be a separable metric space with metric $d_X$, and let $\cS_X=\{s_i: i\in \bN\}\subset X$ be a countable dense subset. We say $(X,d_X,\cS_X)$ is a \textit{computable metric space} if the distance function $d_X(.,.)$ is uniformly computable on $\cS_X\times \cS_X$, that is, if there exists a Turing machine $T=T(i,j,n)$, which on input $i,j,n\in \bN$ outputs a rational number such that
$|d_X(s_i,s_j)-T(i,j,n)|<2^{-n}$.
\end{defn}
The points in $\cS_X$ in Definition \ref{def:computable} are called the \textit{ideal points} of $X$ and $\cS_X$ is the ideal set of the computable metric space. The ideal points assume the role of $\bQ^k$ in $\bR^k$. We may suppress the subscript $X$ and write $(X,d,\cS)$ instead of $(X,d_X,\cS_X)$ when the context precludes ambiguity.

\begin{defn}
Let $(X,d_X,\cS_X)$ be a computable metric space.
An {\em oracle} for $x\in X$ is a function $\phi$ such that on input $n\in \bN$, the output $\phi(n)$ is a natural number so that $d_X(x,s_{\phi(n)})<2^{-n}$.  Moreover, we say $x$ is {\em computable} if there is a Turing machine $T=T(n)$ which is an oracle for $x$.
\end{defn}

Next, we see that Euclidean spaces are computable metric spaces.

\begin{exmp}
Consider the triple $(\mathbb{R}^k,d_{\mathbb{R}^k},\cS_{\mathbb{R}^k})$ with $d_{\mathbb{R}^k}$ the Euclidean distance on $\bR^k$ and $\cS_{\mathbb{R}^k}=\mathbb{Q}^k$. Then $(\mathbb{R}^k,d_{\mathbb{R}^k},\cS_{\mathbb{R}^k})$ is a computable metric space. In particular, for a real number $x$, an oracle for $x$ is a function $\phi$ such that on input $n$, $\phi(n)$ is a rational number so that $|x-\phi(n)|<2^{-n}$.
\end{exmp}

We extend the notion of computable functions to functions between computable metric spaces as follows.

\begin{defn}\label{def:computablefunction}
Let $(X,d_X,\cS_X)$ and $(Y,d_Y,\cS_Y)$ be computable metric spaces and  $\cS_Y=\{t_i: i\in \bN\}$.  Let $D\subset X$.  A function $f:D\rightarrow Y$ is {\em computable} if there is a Turing machine $T$ such that for any  $x\in D$ and any oracle $\phi$ of $x$, the output $T(\phi,n)$ is a natural number satisfying $d_Y(t_{T(\phi,n)},f(x))< 2^{-n}$.

\end{defn}

For example, if $D\subset\mathbb{R}^\ell$ and $f=(f_1,\dots,f_k):D\rightarrow\mathbb{R}^k$, then $f$ is computable if and only if all functions $f_i$ are computable.  We observe that, in this definition, $x$ does not need to be computable, i.e., the oracle $\phi$ does not need to be a Turing machine.  In the case where $x$ is computable, however, $f(x)$ is computable because $T(\phi,n)$ is an oracle Turing machine for $f(x)$.  

The composition of computable functions is computable because the output of one Turing machine can be used as the input approximation for subsequent machines.  In addition, basic operations, such as the arithmetic operations and the minimum and maximum functions, are computable. We refer to \cite{BHW2008} for more details on these topics.

Since the definition for a computable function uses any oracle for $x$ and applies even when $x$ is not computable, we can conclude that for any sufficiently close approximation $y$ to $x$, $f(y)$ approximates the value of $f(x)$; that is, $f$ is continuous. Thus, for a function to be computable, it must be continuous.

\begin{lem}[{\cite[Theorem 1.5]{BY2009}}]\label{lem:continuous}
Let $(X,d_X,\cS_X)$ and $(Y,d_Y,\cS_Y)$ be computable metric spaces, $D\subset X$, and $f:D\rightarrow Y$.  If $f$ is computable, then $f$ is continuous.
\end{lem}

%-------------------------------------------------
\subsection*{Computability of sets, and the $L^\infty$ metric}

Next, we introduce the computability of compact subsets of computable metric spaces.
Recall that the {\em Hausdorff distance} between two {compact} subsets $A$ and $B$ of a metric space $X$ is given by
$$
d_H(A,B)=\max\left\{\max_{a\in A}d(a,B),\max_{b\in B}d(b,A)\right\},
$$
where $d(x,C)=\min\{d_X(x,y): y\in C\}$.   In other words, the Hausdorff distance is the largest distance of a point in one set to the other set. 

\begin{defn} \label{defn:2^n-approx}
    A set $\Psi_n$ is called a \textit{$2^n$-approximation of $C$} if, in the Hausdorff metric, $d_H(C,\Psi_n)\leq 2^{-n}$. 
\end{defn}

\noindent \textbf{Notation.}  
We use the notation $B(s,\delta)$ for the $\delta$-ball about a point $s$, and the notation $\cN(S,\delta)$
for the $\delta$-neighborhood about a set $S$.

\begin{defn}
\label{defn:S_Y}
Let $(X,d_X,\cS_X)$ be a computable metric space.
We say that a ball $B(x,r)$ is an \textit{ideal ball} if $x\in \cS_X$ and $r=2^{-i}$ for some $i \in \bZ$.

 For 
    $\sC= \{C\subset X\,\, {\rm compact}\}$ let $\cS_\sC = \cS_\sC(X)$ denote the collection of finite unions of closed ideal balls.
\end{defn}

We have the following, see e.g., \cite{BY2009}:

\begin{lem}\label{lem:compact-subsets-computable-metric-space}
Let $(X,d_X,\cS_X)$ be a computable metric space, and let $\sC$ and $\cS_\sC$ be as in Definition~\ref{defn:S_Y}. 
Then 
$(\sC,d_H,S_\sC)$ is a computable metric space.
\end{lem}

While considering the computability of sets in $\bC^\ell$ we slightly deviate from \cite{DY2018}  where the Euclidean metric is used. This is because it  is more natural for computer calculations in $\bC^\ell$
to consider vectors in $\mathbb{R}^{2\ell}$
rather than $\CC^\ell$, and use the $L^{\infty}$ metric, rather than 
Euclidean.

\medskip

\noindent \textbf{Notation.}  
When we write $\norm{\cdot}$ we mean the $L^{\infty}$ 
norm on~$\mathbb{R}^{2\ell}$, so that for a vector $z=(z_1,\ldots,z_\ell)\in \mathbb{C}^\ell,$
\begin{equation}
\label{eqn:boxnorm}
\snorm{z}=\max\{ \abs{\Re(z_j)},\abs{\Im(z_j)} \colon 1\leq j\leq \ell\}.    
\end{equation}
Hence, 
$$d(z,w) = \max \{ | \Re(z_j)-\Re(w_j)|, |\Im(z_j)-\Im(w_j)|: 1\leq j \leq \ell \}, $$ 
and if $z$ and $w$ are in a (closed) box of sidelength $r$, then $d(z,w) \leq r$. 
We may use the simpler notation $\sonorm{\cdot}$ in one complex dimension.

The $L^\infty$ metric is {\em uniformly
equivalent} to the euclidean metric on~$\mathbb{C}^\ell, \enorm{\cdot}$, since
$
\frac{1}{\sqrt{2\ell}} \enorm{x} \leq \snorm{x} \leq \enorm{x}.
$
Neighborhoods are slightly different concerning two uniformly equivalent
norms, but the topology generated by them is the same; thus, they
can practically be used interchangeably.

We may say \textit{box} when we mean a ball around a point in the $L^\infty$ norm.  Analogously to \cite{DY2018}, we consider ideal balls to have dyadic rational side length and dyadic rational center coordinates. Indeed,
for $\ell\in \mathbb{Z}^+,$ we consider the computable metric space $(\bC^\ell,d,\cS)$, where $d$ is the $L^\infty$ metric
and $\cS$ is the set of points $z=(z_1,\ldots,z_\ell)$ with the real and imaginary parts of each of $z_i$ dyadic rationals; i.e., points in $\cD = \{ a/2^b : a\in\bZ, b\in \bN\}$.

From now on, we focus primarily on $\CTwo$ as a function space, although we consider a higher-dimensional $\CC^\ell$ as a parameter space. 

Similarly to \cite{Braverman-Yampolsky-2008} aside from using the $L^\infty$ metric,
applying Definition~\ref{defn:S_Y} and Lemma~\ref{lem:compact-subsets-computable-metric-space} to $\bC^2$ we conclude:

\begin{cor} \label{cor:set-computable}
A compact set $C \subset \bC^2$ is \text{computable} if there is a Turing Machine $T(n)$ which on input $n \in \bN$ outputs a set $\Psi_n$, defined as the union of a finite collection of ideal centers and ideal radii forming balls,  which satisfies $d_H(C,\Psi_n)\leq 2^{-n};$ i.e., 
$\Psi_n\in \cS_\sC$ is a \text{$2^n$-approximation of $C$}.    
\end{cor}

Again, in our setting, the Hausdorff metric is based on the $L^\infty$ metric. This makes sense because of how a computer actually draws a picture of a set. Round pixels have to overlap in order to form a cover of Euclidean space. Computer-implemented pixels are squares (not round), and form a perfect, uniformly-sized grid or lattice, only overlapping on their boundaries.

\smallskip

Finally, to study the locus of maps for which an interesting invariant set has some stable behavior (in our case, the locus of Axiom A polynomial skew products of a fixed degree), we use the following. 

\begin{defn}
\label{defn:LowerComputableSet}
Let $D\subset \CC^{\ell}$ be open. We say $D$ is {\em lower semi-computable} if there exists a Turing machine producing $T =\{(x_i,r_i)\}_{i\in\bN}$ such that $x_i \in \CC^\ell$ has dyadic rational coordinates, $r_i$ is a non-negative dyadic rational radius, and $D = \cup_{i=1}^\infty B(x_i,r_i)$. 
\end{defn}

Note that in the definition of lower semi-computability, we do not require an error estimate for how close any finite union of dyadic balls is to the set $D$, just that it converges in the limit. 

%-----------------------------------
\subsection{Invariant sets and hyperbolicity}
\label{sec:dpspaces}
%-----------------------------------
We briefly review some key concepts in dynamical systems. First, we define certain invariant sets that capture the complicated dynamics of the system. 
The {non-wandering set}, $\Omega$,
the {chain recurrent set}, $\cR$, and
the {Julia set}, $J$,   
are ways to
identify (in a certain sense) the set of points with dynamically interesting behavior.

\begin{defn} \label{defn:invariant_sets}
Let  $f:\bC^2\to \bC^2$ be a polynomial endomorphism.

(1) The \textit{non-wandering set} $\Omega=\Omega_f$ is the set of points $x$ such that for every open set $U$ containing $x$, and every $N>0$, there is an $n>N$ such that $f^n(U)\cap U \neq \emptyset$.

The compact,  invariant, and transitive subsets of $\Omega$ are called the  {\em basic sets} of $f$.

(2) The {\em Julia set}, $J$, of $f$ is the topological boundary of the filled-in Julia set, $K$, the set of points in $\CTwo$ with bounded orbits under $f$.

(3) An {\em $\ep$-chain} of length $n>1$ from $y$ to $z$ is a 
  sequence of points
  $\{y=x_1, \ldots , x_n=z\}$ such that  $\abs{f(x_k) - 
  x_{k+1}} < \epsilon$ for $1 \leq k \leq n-1.$
A point $y$ belongs to the {\em $\ep$-chain recurrent set}, 
  $\chreceps$, of $f$ if there is an
  $\ep$-chain from $y$ to $y$.
The {\em chain recurrent set} of $f$ is $\cR = \cap_{\ep>0} 
\chreceps.$

A point $z$ is in the {\em forward chain limit set of a point $y$},
$\cR(y)$, if for all $\ep >0$, for all 
$n\geq 1$, there is an $\ep$-chain from $y$ to $z$ of length 
greater than $n$.
Define an equivalence relation on $\cR$ by: 
$y \sim z$ if $y \in \cR(z)$ and $z\in\cR(y)$.
The corresponding equivalence classes are {\em chain transitive components}, or simply {\em chain components}. 
Analogously, define $\chreceps(y)$ and $\ep$-chain
components. 
\end{defn}

The chain recurrent set $\cR$ is closed and $f$-invariant. If $0<\ep < \ep'$,
then $\cR \subset \chreceps
 \subset \chrecepsp$.  Moreover, $J$ and $K$ are compact,
and $J \subset \mathcal{R}$.

Note that an $\ep$-chain recurrent set depends on the choice of norm, but in any uniformly equivalent norms, the chain recurrent set is the same. 
 Consider  $\cR^e = \cap_{\ep>0} \cR^e_\ep$ defined by the Euclidean metric, and we have $\cR = \cap_{\ep>0} \cR_\ep$ defined by the $L^\infty$ metric on $\CTwo$ viewed as $\bR^4$. Because $\cR^{e}_{\ep} \subset 
\cR_{\ep} \subset \cR^{e}_{\sqrt{2\ell}\ep}$ on $\CC^\ell$, 
the intersections of these nested sequences are the same:
 $\cR = \cR^e$.

 One can sketch an approximation of the intersection of the Julia set  with a complex line by a computer program (testing finite orbits against an escape radius), but
$\cR$ is the set most amenable to
rigorous computer investigations. 
$\cR$ can also be easily decomposed into its chain components, which do not interact with each other. 
Since $\cR \supset J$, 
we can deduce information about $J$ by studying $\cR$.
Finally we note, the set of periodic points of $f$, denoted by Per$(f)$, satisfies $\overline{\text{Per}(f)} \subset \Omega_f \subset \cR_f$.

\subsection*{Hyperbolicity.}
 A map of several variables is {\em hyperbolic} on a compact invariant set $X$ if there exists a continuous splitting of the tangent bundle over $X$ into two subspaces (of any dimension including zero), with one subspace uniformly expanded by the map, and the other uniformly contracted. 
More precisely:

\begin{defn} \label{defn:hyp}
Let $g$ be a $C^1$ diffeomorphism or endomorphism of a compact manifold $M$, and let $X$ be a compact $g$-invariant set. 
Say $X$
is a {\em hyperbolic} set of $g$ if there is a splitting of the tangent bundle $T_xM
= E^s_x \oplus E^u_x$ (one subspace may be trivial), for each $x$ in $X$, which varies
continuously with $x$ in $X$, constants $c>0$ and $\lambda > 1$, and a
Riemannian metric $\norm{\cdot}$ such that the following holds:
\begin{enumerate}
 \item The splitting is $Dg$ invariant: $D_xg(E^s_x) = 
E^s_{g(x)}$, and $D_xg(E^u_x) = E^u_{g(x)}$, and 
\item The tangent map $Dg$ expands (contracts) $E^u (E^s)$ uniformly, that is,
$\norm{D_xg^n(\mathbf{w})} \geq c \lambda^n 
\norm{\mathbf{w}}$ for all $\mathbf{w} \in E^u_x$, and 
 $\norm{D_xg^n(\mathbf{v})} \leq c^{-1}\lambda^{-n} \norm{\mathbf{v}}$ for all $\mathbf{v} \in E^s_x$, for all $n\in\bN$.
\end{enumerate}
We say that $g$ is {\em Axiom A} if the non-wandering set $\Omega_g$ is a hyperbolic set of $g$
and $\overline{\text{Per}(g)}=\Omega_g$. 

\end{defn}

We note that this definition is independent of the choice of Riemannian norm. It is well-known that there always exists a metric for which one can choose $c=1$. Such a metric is called an adapted metric.

A polynomial map in one complex dimension is called {\em hyperbolic} if it is uniformly expanding on its Julia set. 
For a hyperbolic complex polynomial $f$, the Julia set $J$ is a chain component of $\mathcal{R}$, and the other chain components are the attracting cycles of $f$. 

\medskip

\textbf{Notation.} Finally, we note that we use the \textit{infimum} norm $\inorm{A}$ respectively \textit{supremum} norm $\supnorm{A}$ of a Jacobian matrix, to measure the minimum respectively maximum expansion on tangent vectors $\mathbf{v}$, by the matrix:
$$
\inorm{A} = \displaystyle \min_{\mathbf{v}, \snorm{\mathbf{v}} = 1} \norm{A\mathbf{v}} 
\ \ \text{ and } \ \ 
\supnorm{A} = \displaystyle \max_{\mathbf{v}, \snorm{\mathbf{v}} = 1} \norm{A\mathbf{v}}.
$$
We note the infimum norm is actually not a norm since $\inorm{A}=0$ does not imply $A=0$.
Picture the image of a small ball under a matrix as a small ellipse. The radius of the largest ball that fits in this ellipse is the inf norm of the derivative (times the original ball size), and the radius of the smallest ball containing the image is the sup norm of the derivative (times the original ball size). 
We recall that $\inorm{A} = 1/\supnorm{A^{-1}}$ holds for any invertible matrix $A$.

%----------------------------------------------
\subsection{Polynomial skew products of $\CC^2$}
\label{sec:skewintro}
%-----------------------------------

In this subsection, we first summarize some notation and results (primarily from \cite{Jonsson1999}), to give some background on polynomial skew products: $f(z,w) = (p(z),q(z,w))$, where $p$ and $q$ are both of degree $d\geq 2$.

Since the dynamics of $f$ in the $z$-coordinate is given by $p$, it is useful to employ the notation ($K_p$ and) $J_p$ for the one-dimensional (filled) Julia set of $p$, and $G_p(z)$ for the Green function in $\CC$ of $p$, where $K_p = \{ G_p=0\}$.

\textbf{Global Dynamics.} For polynomial skew products,  the usual rate of escape Green function, defined for $x\in\CTwo$ by $G(x)=\lim_{n\to\infty} \frac{1}{d^n} \log^+ \abs{f(x)}$,  is continuous, plurisubharmonic, nonnegative, and satisfies $G\circ f = dG$ and  $K=\{G=0\}$.
One can also define a positive closed current $T = \frac{1}{2\pi}dd^c G$ and an ergodic invariant measure, $\mu = T \wedge T$, of maximal entropy $\log d^2$.  
 Then $J_2:=\supp(\mu)$ coincides with the closure of the set of repelling periodic points.

\textbf{Vertical dynamics.}  Since $f$ preserves the vertical lines $\{z\} \times \CC$, it is useful to consider the dynamics of $f$ on this family of lines.  Let $z_n = p^n(z)$, $q_z(w) = q(z,w)$,  and $Q^n_z(w) =  \circ \cdots \circ q_z$, so that $f^n(z,w) = (z_n, Q^n_z(w)).$  Let $G_z(w) = G(z,w) - G_p(z)$.  Then $G_z$ is nonnegative, continuous, subharmonic, and is asymptotic to $\log\abs{w} - G_p(z)$ as $w\to\infty$.  
Naturally, define $K_z = \{ G_z=0\}$, and $J_z=\partial K_z$.  Then $K_z$ and $J_z$ are compact, and  if $z\in K_p$, then $w\in K_z$ if and only if $\abs{Q^n_z(w)}$ is bounded.  Further, $G_{z_1} \circ q_z = d G_z$, which implies $q_z(K_z) = K_{z_1}$ and  $q_z(J_z) = J_{z_1}$.

However, not every phenomenon of one-dimensional dynamics carries over to vertical dynamics.  For example, unlike in one dimension, $J_z$ may have finitely many (but greater than one) connected components, even for $d=2$ (see \cite{Jonsson1999}, remark 2.5).

\textbf{Vertical Expansion.}  Let $Z \subset K_p$ be compact with $pZ \subset Z$, for example $Z=J_p$ or $Z=A_p$, the set of attracting periodic orbits.  Let $J_Z = \overline{\cup_{z\in Z}\{z\} \times J_z }$.  Jonsson showed
$J_{J_p} = J_2$.
Call $f$  {\em vertically expanding over $Z$} if there exist $c>0$ and $\lambda >1$ such that $\abs{DQ^n_z(w)} \geq c\lambda^n$, for all $z\in Z$, $w\in J_z$, and $n\geq 1$.

As in one dimension, $z\mapsto K_z$ is upper semi-continuous and $z\mapsto J_z$ is lower semi-continuous, in the Hausdorff metric.  Further, if $f$ is vertically expanding over $Z$, then $z\mapsto J_z$ is continuous for all $z\in Z$; and if in addition, $J_z$ is connected for all $z\in Z$, then $z\mapsto K_z$ is continuous for all $z\in Z$.  However, Jonsson provides examples showing vertical expansion over $J_p \cup A_p$ neither implies that $z\mapsto J_z$ is continuous on all of $\CC$, nor that $z\mapsto K_z$ is continuous for $z\in J_p$.
Jonsson obtained a very useful equivalent condition for $f$ being an Axiom~A polynomial skew product, namely:

\begin{thm}[\cite{Jonsson1999}, Theorem 8.2]
\label{thm:AxiomAVertExp}
A polynomial skew product $f$ is Axiom A on $\CTwo$ if and only if 
\begin{enumerate}
\item $p$ is uniformly expanding on $J_p$,
\item $f$ is vertically expanding over $J_p$, and
\item $f$ is vertically expanding over $A_p$.
\end{enumerate}
Moreover, if $f$ is Axiom A, then $\cR = \Omega_f = \overline{{\rm Per}(f)}$.
\end{thm} 

To reiterate, if $f$ is an Axiom A polynomial skew product, then the non-wandering set equals the chain recurrent set, the basic sets are the chain components, and there is a continuous splitting on the tangent bundle over that set into directions which are uniformly expanded or contracted. Moreover, due the preservation of the vertical fibers, the splitting is restricted, allowing only for uniform expansion in both directions (on one of the chain components, $J_2$), uniform contraction on both directions (on any globally attracting periodic cycles, each is a chain component), or saddle behavior split along horizontal versus vertical (creating two more types of chain components). 

Jonsson also provides a structural stability result for Axiom A skew products (\cite{Jonsson1999}, Theorem A.6 and Proposition A.7).  It follows that being Axiom A is an open condition, and it makes sense to refer to a connected component of the subset of Axiom A mappings in a given parameter space as a {\em hyperbolic component}. 

%----------------------------
See Appendix~\ref{sec:App_Examples} for a description of some different types of Axiom A polynomial skew products. 

%=============================================
\subsection{Preliminaries for computability of skew products}
\label{sec:imagebounds}
Let $d\geq 2$. Recall from Section~\ref{sec:intro-results} that $\mathscr{A}_d$ denotes the parameter space of Axiom A polynomial skew products of $\bC^2$ of degree $d$. 
It is our goal to establish the computability of the following functions: $\mathscr{A}_d\ni c\mapsto J_2=J_2(c)$, $\mathscr{A}_d\ni c\mapsto X_0=X_0(c)$, $\mathscr{A}_d\ni c\mapsto X_1=X_1(c)$, and $\mathscr{A}_d\ni c\mapsto \Lambda=\Lambda(c)$.
These four computability results are evidently stronger statements than Theorem \ref{thm:main} since they include that one can use the same Turing machine to compute $J_2(c), X_0(c),X_1(c)$ and $\Lambda(c)$ independently of $c\in \mathscr{A}_d$. 

We start by establishing some preliminary computability aspects of polynomial skew products. Since $f(z,w) = (p(z),q(z,w))$, where $p$ and $q$ are both polynomials of degree $d$,
we have 
 \begin{equation}
 \label{eqn:coefficients_f}
p(z) = \displaystyle\sum_{0 \leq k \leq d} a_k z^k\  \text{ and  } \
q(z,w) = \sum_{0\leq j+k \leq d} b_{k,j} z^k w^{j},  
\end{equation}
with at least one $b_{j,k}\not=0$ with $j+k=d$.
 We can always assume that $p$ and $q$ are monic, so $a_d = b_{0,d} = 1$, and after a linear change of coordinates, we may assume that $a_{d-1}= b_{1,d-1}=0$ (see \cite{Jonsson1999}, Definition 1.1 and right above Definition 7.1). 
Hence the parameter space of polynomial skew products is $\CC^{\ell}$ where $\ell = (d-1) + \frac{1}{2}(d+1)(d+2) - 2 $.
One basic ingredient for obtaining the computability results for the sets $J_2, X_0, X_1,\Lambda$ is to estimate how close or far an image of a point is from an approximate point. Below, we provide two results that study this for polynomial skew products.  

The following will be  used in Section~\ref{sec:results1} to obtain an error estimate for calculating images of boxes under a polynomial skew product.

\begin{lem}
\label{lem:linearization-error}
    Let $f=f_c$ be a polynomial skew product of $\CTwo$ of degree $d\geq 2$ which depends on a parameter $c$. There exists a Turing machine $T=T(c,z_0,w_0,r)$, which on input of $c$, $(z_0,w_0)\in \CC^2\cap \bQ^4$ and $r\in \bQ^+$ outputs a rational number $L>0$ such that
    \begin{equation}\label{Lipconst1}
        \snorm{f(z_1,w_1)-f(z_2,w_2)} \leq L\snorm{(z_1,w_1)-(z_2,w_2)}
     \end{equation}
    for all $(z_1,w_1), (z_2,w_2)\in \overline{B((z_0,w_0),r)}$.
\end{lem}

\begin{proof}
Recall that $f(z,w)=(p(z),q(z,w))$ where $p$ and $q$ are degree $d$ polynomials depending on the parameter $c$. Since
\[D_{(z,w)}f = \twovec{p'(z) \hskip14mm  0 }{q_1(z,w) \ \ \ \ \ q_2(z,w)}
\] 
where $q_1 = \partial q / \partial z,$ and  $q_2 = \partial q/\partial w,$ it is straight-forward to compute $L\in \bQ^+$ such that  $\supnorm{D_{(u,v)}f}\leq L$ for all $(u,v)\in B((z,w),r)$. Therefore, Equation~\eqref{Lipconst1} follows from the Mean Value Theorem for vector-valued functions.
\end{proof}

We can apply this result to calculate a uniform $L$ which holds for any finite collection of closed balls.

\begin{cor}\label{cor:box-error}
      Let $f=f_c$ be a polynomial skew product of $\CTwo$ of degree $d\geq 2$ which depends on a parameter $c$.
      Let $\cB=B_1\cup\dots \cup B_s$ be a finite union of ideal closed balls in $\CTwo$.
      Then there exists a Turing machine $T=T(c,B_1, \dots, B_s)$, which outputs a rational number $L>0$ such that
    \begin{equation}\label{Lipconst2}
        \snorm{f(z_1,w_1)-f(z_2,w_2)} \leq L\snorm{(z_1,w_1)-(z_2,w_2)}
     \end{equation}
    for any two points $(z_1,w_1), (z_2,w_2)$ which lie together in the same ball $B_i$ in the collection $\cB$.
\end{cor}

\begin{proof} 
In the proof of the above lemma, $L$ is an upper bound for the sup norm of the derivative of $f$ in a ball of interest. Since there are a finite number of balls here, we take the largest $L$ for this corollary. 
\end{proof}

Next, we provide a computational criterion for $f$ being distance expanding, which is utilized in Section~\ref{sec:proof-thm-1}.  

\begin{lem}
\label{lem:distance-expansion_cw}
     Let $f=f_c$ be a polynomial skew product of $\CTwo$ of degree $d\geq 2$ depending on a parameter $c$. Let $(z_0,w_0)\in \bC^2$, $r\in \bQ^+$ and $\lambda>0$  such that $\inorm{D_{(z,w)}f}\geq 1+\lambda$ for all $(z,w)\in \overline{B((z_0,w_0),r)}$. There exists a Turing machine $T=T(c,z_0,w_0,r)$ which on input of $c$, $(z_0,w_0)\in \CC^2\cap \bQ^4$ and $r$, outputs a rational number $r'>0$ such that
     \[
     \norm{f(z_1,w_1)-f(z_2,w_2)}\geq (1+\lambda) \norm{(z_1,w_1)-(z_2,w_2)}
     \]
     for all $(z_1,w_1),(z_2,w_2)\in \overline{B((z_0,w_0),r')}$.
     \end{lem}
     \begin{proof}
     The inverse function theorem guarantees the existence of $r'>0$ such that $f\vert_{\overline{B((z_0,w_0),r')}}$ is invertible with a holomorphic (polynomial) inverse. Once we can show that such a radius $r'$ is computable from the input data, the proof of the lemma is identical to the proof of Theorem 6.4.1 in \cite{Urbanski_book}.  We write $A=D_{(z_0,w_0)}f$. We recall from the proof of the inverse function theorem that $f\vert_{\overline{B((z_0,w_0),r')}}$ is invertible whenever
     \begin{equation}\label{inversefunc}\snorm{A-D_{(z,w)}f}<\frac{1}{2\snorm{A^{-1}}}
  \end{equation}
for all $(z,w)\in \overline{B((z_0,w_0),r')}$. Since $A$ is computable from the input data, it follows that $\snorm{A^{-1}}$ is also computable. We recall that 
 \begin{equation}\label{eqdf}
 D_{(z,w)}f = \twovec{p'(z) \hskip14mm  0 }{q_1(z,w) \ \ \ \ \ q_2(z,w)}
\end{equation} 
Since $p',q_1,q_2$ are polynomials depending on the parameter $c$ it follows that $Df$ has a computable modulus of continuity. Therefore, we can compute $r'>0$ satisfying Equation~\eqref{inversefunc}. This completes the proof.
\end{proof}

By slightly modifying the proof of Lemma \ref{lem:distance-expansion_cw}, we obtain a uniform estimate for the radius $r'$.

\begin{cor}
\label{cor:distance-expansion-cw}
 Let $f=f_c$ be a polynomial skew product of $\CTwo$ of degree $d\geq 2$ depending on a parameter $c$. Let $\cB=B_1\cup\dots \cup B_s$ be a finite union of ideal closed balls and $\lambda>0$  such that $\inorm{D_{(z,w)}f}\geq 1+\lambda$ for all $(z,w)\in \cB$. Then there exists a Turing machine $T=T(c,B_1,\dots,B_s)$ which outputs a rational number $r'>0$ such that
     \[
     \norm{f(z_1,w_1)-f(z_2,w_2)}\geq (1+\lambda) \norm{(z_1,w_1)-(z_2,w_2)}
     \]
     for all $(z_1,w_1),(z_2,w_2)\in B_i$ for  $i=1,\dots,s$ with $\norm{(z_1,w_1)-(z_2,w_2)}\leq 2r'$.   
\end{cor}
\begin{proof}
The proof is analogous to the proof of Lemma \ref{lem:distance-expansion_cw}. The only required modification is to make sure that Equation~\eqref{inversefunc} holds for all $(z,w)\in  \overline{B((z_0,w_0),r')}$ independently of $(z_0,w_0)\in \cB$. This follows from the following argument: Since  $(z,w)\mapsto D_{(z,w)}f$ has a computable modulus  of continuity on $\cB$, the map $(z,w)\mapsto \snorm{(D_{(z,w)}f)^{-1}}$ also has a computable modulus of continuity on $\cB$. Thus, we can compute $\gamma>0$ such that  $\snorm{(D_{(z,w)}f)^{-1}}\leq \gamma$ for all $(z,w)\in \cB$.
Therefore, we can consider a modified version of Equation~\eqref{inversefunc}. Namely, it suffices to show that one can compute $r'>0$ such that
\begin{equation}\label{inversefunc1}\snorm{D_{(z_0,w_0)}f-D_{(z,w)}f}<\frac{1}{2\gamma}.
  \end{equation}
 holds for all  $(z_0,w_0)\in \cB$ and all $(z,w)\in \overline{B((z_0,w_0),r')}$, which follows from the computable modulus of continuity of the map $(z,w)\mapsto \snorm{D_{(z,w)}f}$.
 
\end{proof}

Above, we used $2r'$ (instead of $r'$) for consistency with \cite{Urbanski_book}, Definition 4.1.1.

%====================================
\section{Computability of $\mathcal{R}$ for polynomial skew products}
\label{sec:results1}
%====================================

%--------------------------------------
\subsection{Computing $\cR_\alpha$ for polynomial skew products}
\label{sec:hypatia}
%--------------------------------------

In \cite{Hruska-Skew-2006}, the first author of this paper (under a former name) describes a rigorous computer algorithm and its implementation for polynomial skew products of $\CTwo$, constructing a decreasing sequence of nested neighborhoods $\cB_n$ of 
the chain recurrent set, $\cR$,
and a graph $\Gamma_n$ modeling the dynamics of $f$ on $\cB_n$. 

While we don't use this aspect of the work for our computability results, \cite{Hruska-Skew-2006} contains an algorithm to rigorously test for Axiom A, in the sense that if the computer program determines a map is Axiom A, then it is. This result extends techniques of previous works of the same author \cite{Hruska-HenonHyp,Hruska-HenonJ}.  
 
A similar approach was previously applied in different settings
to develop a very general procedure for rigorously approximating $\chrec$ for continuous maps or flows in $\mathbb{R}^k$, including: Osipenko and Campbell (\cite{Osi, OsiCamp}) approximate the chain recurrent set for a homeomorphism of a smooth, real, compact manifold, Eidenschink (\cite{Eiden}) discusses a similar procedure for real flows, and a philosophically related procedure is studied in \cite{DellHoh, Dell1}, though their case of interest is the attractor of a real map (rather than the chain recurrent set).  The article \cite{KMisch} surveys results in this direction prior to \cite{Hruska-HenonJ}.

The basic idea in these works is to construct a graph that is called a {\em box chain recurrent
model}, which we define below. It satisfies Osipenko's definition of a {\em symbolic image} of $f$ (\cite{Osi}). 

%-----------------

This algorithm to compute a neighborhood of $\cR$ produces the following. 

\begin{defn} \label{defn:boxchrecmodel}
Let $\cR$ be the chain recurrent set of a map $g \colon \mathbb{C}^\ell \to \mathbb{C}^\ell$.
Let $\Gamma = (\mathcal{V}, \mathcal{E})$ be a directed graph,
with vertex set $\mathcal{V} = \{ B_k \}_{k=1}^m$,  where the sets $B_k$ are closed boxes in $\mathbb{C}^\ell$ with pairwise disjoint interior, and such that the union of the boxes $\cB = \cup_{k=1}^m B_k$ contains $\cR$.   
We further require that there is an edge from $B_k$ to $B_j$ if 
$g(B_k)$ 
intersects $B_j$, that is,  
\[
 \{ (k, j) \colon 
g(B_k)
\cap B_j  \neq \emptyset \} \subset \mathcal{E}(\Gamma).
\]
We also require that  $\Gamma$ is the disjoint union of strongly connected components $\Gamma'_i$, $i=1,\dots,s$; that is, for $i\in 1,\dots,s$ and any $B_k, B_j\in \Gamma'_i$, there is a path in $\Gamma'_i$ from $B_k$ to $B_j$.  

If these properties hold
we say that $\Gamma$ is a {\em box chain recurrent model} of $g$ on $\cR$, and the $\Gamma'_i$ are the {\em box chain components} of the model.
\end{defn}

A box chain recurrent model $\Gamma$ provides an approximation of the dynamics of $g$ on $\cR$, 
and the connected components $\Gamma'_i$ are approximations of the chain components.

Here, we use the same general outline as in \cite{Hruska-Skew-2006} for the (optionally, inductive) process for creating $\Gamma$ for polynomial skew products, except we have to modify the approach to make it suitable to establish computability. Precisely, we have to ensure that whenever the input map is Axiom A, the subsequent algorithms which utilize this construction will halt, establish hyperbolicity, and output approximations of $J_2,X_0,X_1,\Lambda$ (whichever exist) at any given pre-described accuracy. Since we are not implementing this algorithm, we also perform some simplifications to ease the exposition.

\begin{alg}
\textbf{Constructing a box chain model,  of a given level $n$, for a polynomial skew product.} 
\label{alg:box-cover}
Let $f(z,w) = (p(z),q(z,w))$ be a polynomial skew product, with $p$ and $q$ of the same degree $d\geq 2.$  For any positive integer $n$, we construct a graph $\Gamma_n$ satisfying Definition~\ref{defn:boxchrecmodel}, whose vertices are a finite collection of (closed) boxes $\cB_n = \cup_l B^n_k$. 

Additionally, our $\cB_n$ has the property that for a user-inputted $\xi_n>0$, it is guaranteed that there is no edge $(k,j)$ from the box $B_k^n$ to the box $B_j^n$ whenever $d_H(f(B_k),B_j)> \xi_n.$
    
(1) Start by \textbf{calculating $R>0$} such that $\cR \subset [-R, R]^4$. This calculation is similar to one in \cite{Hruska-HenonJ}. Then, for simplification in our computability argument, let $R$ be the smallest power of $2$ which suffices, say $R=2^m$. We start with $\cB_0 = [- R, R]^4.$ 

(2) For any positive integer $n$ desired, evenly partition $\cB_0 = [-R,R]^4$ into a grid of $(2^n)^4$ two-complex-dimensional boxes
$B^n_l = B^{n}_{k,j} = B^{n,z}_k \times B^{n,w}_j$, where the boxes
$B^{n,z}_1,\dots, B^{n,z}_{4^{n}}$  (and, respectively, $B^{n,w}_1,\dots, B^{n,w}_{4^{n}}$) are in the square $[-R, R]^2$ in the $z$-plane (and respectively, $w$-plane).  

Each box $B^{n,z}_k$ and $B^{n,w}_j$, and thus each $B^n_l$, has sidelength (thus diameter in our $L^\infty$ norm)  a dyadic rational and center point a dyadic rational. Let $\ep_n = 2R/2^{n} = 2^{-n+m+1}$ denote the box sidelength.

(3) Next, calculate the approximate image of each box $B^n_l$ given the prescribed accuracy $\xi_n$. 
 To do so, partition each box $B^n_l = B^{n}_{k,j} = B^{n,z}_k \times B^{n,w}_j$ by splitting each $B^{n,z}_k$ and each $B^{n,w}_j$ (including their boundaries) into a grid of $2^{s_n}\times 2^{s_n}$ ideal (dyadic rational) points, for some (positive) integer $s_n$.
 
Calculate which $s_n$ is needed based on $\xi_n$ by applying Corollary~\ref{cor:box-error} to the collection of all boxes $\cup_l B^n_l$. The corollary guarantees that we can compute a constant $L=L_n$  such that for any two points in the same box, any $B^n_l$, the distance between their images is bounded by a factor $L_n$ times the distance between the points.  

Now if $x'$ and $y'$ are adjacent partition points, $d(x',y') = 2^{-s_n -n+m-1}$. 
So, given $\xi_n>0$, we need
$s_n$ to satisfy:
\\
$d(f(x'),f(y')) \leq L_n \cdot  d(x',y') = L_n \cdot 2^{-s_n -n+m-1} \leq \xi_n$. Since $L_n$ and $\xi_n$ are determined, set $s_n$ to be the smallest integer so that 
$L_n\cdot 2^{-s_n -n+m-1} \leq \xi_n$.

Let $\sG_{\xi_n}(B^n_l) = \{ x' = (z',w') $ in our $(2^{s_n})^4$-grid of $B^n_l  \}, $ which includes partition points along all box boundaries (i.e., if we are partitioning, a real interval $[a,a+\delta]$ by $2^s$ points, they are the points $\{a\cdot (j \delta/(2^s-1)) + (a+\delta) \cdot (1-(j \delta/(2^s-1))), j=0,\ldots s$\}).
Now that $s_n = s_n(\xi_n)$ is determined, we  
define, and then \textbf{calculate} for each box $B^n_l$, the approximate image $F_{\xi_n}(B^n_l)$ as the union of the (closed) boxes (balls in our norm) of radius $\xi_n$ with centers all of the $f(x')$ for each grid point $x'$ in the box $B^n_l$; that is,

\begin{equation}
\label{eqn:F=approx-image-of-f}
F_{\xi_n}(B^n_l) := \bigcup_{x'\in \sG_{\xi_n}(B^n_l) } \overline{B(f(x'),\xi_n)}    
\end{equation}

Observe that 
$d_H(F_{\xi_n}(B^n_l),f(B^n_l)) \leq \xi_n$, which follows from the calculation of $s_n$ based on $\xi_n$, as described above.

(4) With these approximate images of each box calculated, build a transition graph for $f$ which we call $\Upsilon = \Upsilon_n(f)$, whose vertices are the boxes $B^n_l$, and
  \textbf{create an edge $(k,j)$ if the approximate image $F_{\xi_n}(B^n_k)$ lies within $\xi_n$ of $B^n_j$}; i.e., if $F_{\xi_n}(B^{n}_{k}) \cap \cN(B^n_j,{\xi_n}) \neq \emptyset$. Then there is definitely an edge $(k,j)$ if $f(B^{n}_{k}) \cap B^n_j \neq \emptyset$, based on the accuracy of the approximation $F_{\xi_n}$ to $f$ as stated in step (b) above,
and there is no edge from 
 $B_k^n$ to $B_j^n$ if $f(B_k^n) \cap \overline{\cN(B_j^n,{\xi_n})} = \emptyset$, i.e., $d_H(f(B_k^n), B_j^n) > \xi_n$.

(5) Now that the graph $\Upsilon_n$ is formed,  \textbf{compute} the maximal subgraph $\Gamma_n$ of $\Upsilon_n$ which consists precisely of edges and vertices lying in cycles.
This can be done by finding the Strongly Connected Components of the graph, using one of the standard algorithms, see e.g. \cite{CLR}. For example, Kosaraju's algorithm and Tarjan's algorithm both rely on Depth First Search and have a time complexity of $O(\mathcal{V}+\mathcal{E})$, where $\mathcal{V}$ is the number of vertices and $\mathcal{E}$ is the number of edges of $\Upsilon_n$.

The boxes in $\CTwo$ associated with $\Gamma_n$ form the set $\cB_n$, which we call the level $n$-approximation. Its boxes have sidelength (diameter) $\ep_n = R/2^{n-1}=2^{m-n+1}$ and radius $R/2^n=2^{m-n}$, which is again (using $R=2^m$) a power of $2$. Hence, the boxes of the level $n$ approximation have a dyadic rational sidelength, and their center coordinates are also dyadic rationals. 

(6) Note the graph $\Gamma_n$ is partitioned by a disjoint union of $\Gamma_n$'s edge-connected components. 
\textbf{Decompose} $\Gamma_n$ into its edge-connected components; it follows that each of these  
is a box chain component. 

(7) The final boxes are a subset of a $(2^n)^4$-grid on $[-R, R]^4$. If, after all of this, a refined approximation is desired, one may start from the beginning with an increased $n$, or \textbf{iterate} the algorithm by subdividing the boxes at one level $\cB_{n-i}$ to find the boxes at a higher level $\cB_{n}$. Iterating the algorithm by subdividing boxes is more efficient than beginning with a $(2^n)^4$-grid on $[-R, R]^4$ for very large $n$. 

In the algorithms in the next subsection, we frequently apply this algorithm to iteratively refine box collections. 
\qed
\end{alg}

Note that since the finite set (for a fixed $n$) of constructed (closed) boxes have dyadic rational coordinate centers and all the same dyadic rational sidelengths, these are all ideal balls in this norm.

Algorithm~\ref{alg:box-cover} differs from that of \cite{Hruska-Skew-2006} in a couple of ways. First  \cite{Hruska-Skew-2006} used Interval Arithmetic and calculated Hull$(f(B))$, a rectangle containing $f(B)$, and drew an edge anytime the Hull of the image of one box intersected another box. But for our purposes, we cannot allow for such an uncontrolled over-estimate, so we had to develop a more refined algorithm to calculate the images of the boxes. Secondly, since \cite{Hruska-Skew-2006} involved implementation, efficiency was important. For example, a more efficient algorithm would start by calculating $R_1>0$ such that $\cR(p) \subset [-R_1, R_1]^2$, and  $R_2>0$ so that $\mathcal{R}(f) \subset [-R_1, R_1]^2 \times [-R_2, R_2]^2$, then performing the 1-dimensional version of the algorithm in the $z$-plane to get a set of boxes in $[-R_1,R_1]^2$, then building boxes in $\CTwo$ in $[-R_2,R_2]$ in the $w$-plane only over that base set of refined boxes. There is no need for the size of the $z$-grid and $w$-grid to be the same when implementing the algorithm, but it eases the exposition.  

\medskip

Next, we show our algorithm produces points which are all $(\ep'_n+2\xi_n)$-pseudo-periodic in our box norm, for a computable $\ep'_n$ and given $\xi_n$. 
We refer to \cite{Hruska-HenonJ} for a similar result for \Henon maps.

\begin{lem} \label{lem:imageboxsize}
Let $\Gamma_n$ be the box chain model of a degree $d$ polynomial skew product $f$ produced by Algorithm~\ref{alg:box-cover}. 
Then there exists a Turing machine $T = T(n,f,R,\xi_n)$ which outputs a rational number
$L_n>0$ such that for $\ep_n' :=(L_n + 1)\ep_n$,
and for any $B_k \in \cB_n = \mathcal{V}(\Gamma_n)$, we have:
\begin{enumerate}
\item  diam$(F_{\xi_n}(B_k))\leq  L_n  \cdot $diam$(B_k)= L_n \ep_n$,
  and
  
\item if $(k,j)\in\mathcal{E}(\Gamma_n)$ (i.e., there is an edge in $\Gamma_n$ from box $B_k$ to box $B_j$), then
for any $x_k \in B_k$ and any $x_j \in B_j$, 
$d(f(x_k), x_j) \leq \ep'_n+2\xi_n$. 
\end{enumerate}
\end{lem}

\begin{proof}
Using $\ep'_n+2\xi_n = (L_n + 1)\ep_n +2\xi_n =L_n\ep_n + \ep_n +2\xi_n, $
the second item follows from the first item, and the fact
that we defined $\Gamma_n$ so that if
$\cN(B_j,{\xi_n}) \cap {F_{\xi_n}(B_k)} \neq \emptyset$, where $d_H(F_{\xi_n}(B_k),f(B_k))\leq \xi_n$, then there must be an edge from $B_k$ to $B_j$. The $L_n \ep_n$ is from the first item, the extra $+\ep_n$ in the sum is from the width of $B_k$, one $\xi_n$ is from the radius of the balls about the images of the partition points on $B_k$, and the other $\xi_n$ is the distance allowed from $B_j$ to $F_{\xi_n}(B_k)$. 

To prove the first statement, apply Corollary~\ref{cor:box-error} again with $r=\ep_n$, so that any two points in a common box have the distance between their images a multiple $L_n$ of the distance between the points. 
But then applying this to points on opposite boundaries of the box yields the result, setting $L_n$ to be the $L$ guaranteed by the corollary.  

\end{proof}

%------------------------------------
The previous lemma immediately implies the following.

\begin{cor} \label{cor:imageboxsize}
For $\ep_n,\ep'_n, L_1 $ as in Lemma~\ref{lem:imageboxsize},  for every $n>1$, $\ep'_n$ decreases to $0$ as $n\to\infty$, with 
$\ep_n < \ep'_n \leq (L_1+1)\ep_n$.
\end{cor}

\begin{proof}
Note that since $\ep'_n=(L_n+1)\ep_n$ we have $\ep_n < \ep'_n$. 
By Lemma~\ref{lem:linearization-error}, $L_n$ is the largest sup norm of $Df$ on any box in the collection $\cB_n$. Hence $L_n$ does not increase with $n$, since the boxes are nested $\cB_{n-1} \subset \cB_n$ by construction. Thus $L_n \leq L_1$ for all $n$.  Then $\ep'_n = (L_n + 1)\ep_n \leq (L_1+1) \ep_n$ for all $n$.
\end{proof}

Since each box in $\Gamma_n$ must lie in a cycle, the two prior results imply the following.

\begin{cor}
\label{cor:pseudo-orb-calc}
For sufficiently large $n$ and for any $\alpha > \ep'_n+2\xi_n$, all points in $\cB_n$ are $\alpha$-pseudo-periodic, that is,  $\cB_n \subset \cR(\alpha)$. 
\end{cor}

Note that by Corollary~\ref{cor:imageboxsize}, $\ep'_n+2\xi_n$ is computable and tends to $0$ as $n\to\infty$, so long as we choose a computable sequence $\xi_n \to 0$.

%-------------------------------------------

Thus, we have the following result. 

\begin{thm}
\label{thm:boxchainconstruction}
Algorithm~\ref{alg:box-cover} is a Turing machine, depending on $n$ and $f$, and a user-inputted sequence of ideal positive rational $\xi_n \downarrow 0,$ which constructs a finite collection of closed boxes $\mathcal{B}_n =\cup_k B_k^{n} $, 
satisfying:
 \begin{enumerate}
   \item  The boxes $B^n_k \in \cB_n$  are of sidelength $\ep_n=R/2^{n-1}=2^{m-n+1}$, which is dyadic rational. 
   Moreover, by use of the $L^{\infty}$ norm, the diameter of a multi-dimensional box coincides with the sidelength. 

\item  There is an edge in $\Gamma_n$ from $B_k^n$ to $B_j^n$ if $f(B_k^n)) \cap B_j^n \neq \emptyset$, and there no edge from 
 $B_k^n$ to $B_j^n$ if $\overline{\cN(f(B_k^n),{\xi_n})} \cap B_j^n = \emptyset$. 

\item Every box in $\cB_n$ lies in a cycle in the graph $\Gamma_n$.

\item
There is a Turing machine $T=T(n,f,R,\xi_n)$ which outputs
 $ \ep'_n>0$, with $\ep'_n$ non-increasing and tending to zero as $n\to \infty$, such that  
 for any $\alpha > \ep'_n+2\xi_n$, we have $ \mathcal{R} \subset  \mathcal{B}_n \subset \mathcal{R}(\alpha) $.
\\ Hence, $\cR = \cap_{n} \cB_n$. 
 \item The boxes are nested upon refinement: $\mathcal{B}_n \subset \mathcal{B}_{n-1}$ for all $n\geq 0$.
 \end{enumerate}

  \end{thm}

  The $\xi_n \downarrow 0$ chosen by the user should be small compared to the box sidelength $\ep_n$, e.g., $\xi_n = \ep_n/16$ works. 

We close this section with an example $\cB_n$ from an actual implementation for a polynomial skew product, from \cite{Hruska-Skew-2006}. 

\begin{exmp} \label{exmp:hypatia-skew} 
The map $f(z,w) = (z^2-90, w^2 + z/6 + 1.4 + 0.75i)$ satisfies the following: 
\begin{enumerate}
\item $J_p \subset D_1 \cup D_2$, for the intervals $D_1 = [-\beta, -\eta]$ and $D_2 = -D_1$, where $\alpha =-9$ and $\beta=10$ are the fixed points of $p(z)=z^2-90$,
and $\eta = \sqrt{90-\beta}$;
\item for all $z\in D_1$, $q_z$ is in the ``rabbit'' hyperbolic component of the Mandelbrot Set, $\mathcal{M}$, (hence $q_z(w)$ is topologically conjugate to Douady's Rabbit, $\approx w^2 -0.122561+0.744862i$), for which the critical point $0$ lies in a (super)-attracting period three cycle;
and
\item for all $z\in D_2$, $q_z$ is outside of $\mathcal{M}$, i.e., the critical orbit of $q_z$ is unbounded.
\end{enumerate}   

For this map, $J_z$ in the fiber over the $\alpha$-fixed point is the Rabbit, while $J_z$ in the fiber over the $\beta$-fixed point is a Cantor set.   Note that since $A_p$ is empty, the only chain components are the expanding set $J_2$, and a saddle set $\Lambda$ over $J_p$. 
In the fiber over the $\alpha$-fixed point, the attracting 3-cycle of the Rabbit yields a saddle $3$-cycle in $\CTwo$.  This 3-cycle is contained in its own chain component, $\Lambda$, which the algorithm must separate from the component containing $J_2$ before vertical expansion can be established. Our implementation to separate the chain components required
boxes from a $(2^{10})^2 \times (2^9)^2$ grid on $\cB = [-10.1, 10.1]^2 \times [-2.426, 2.426]^2$. 
 Figure~\ref{fig:matrab} shows the components for the box graph model, of size $(|\mathcal{V}|;|\mathcal{E}|)=(78{,}994; 2{,}066{,}558)$, 
which was used to prove that this map is Axiom A as described in \cite{Hruska-Skew-2006} (via a different algorithm than we discuss in the present article).

\begin{figure}
\includegraphics[width=.95\textwidth]{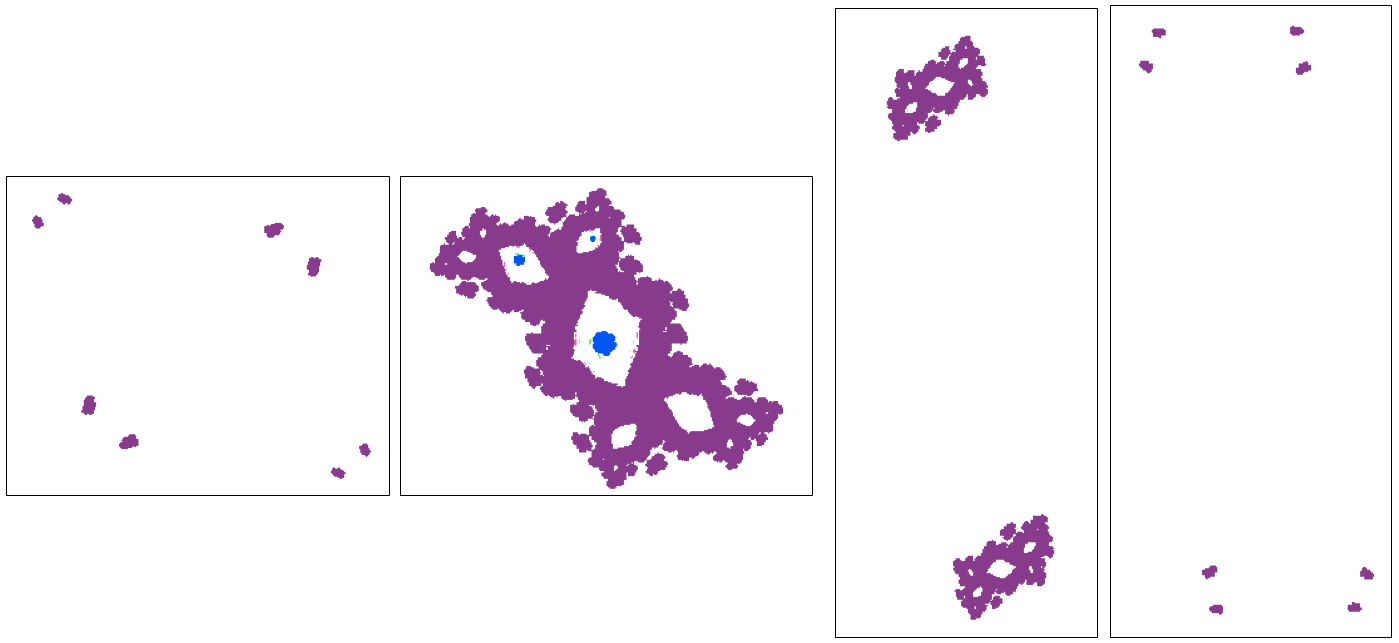}
\caption{For the map $(z,w) \to (z^2-90, w^2+ z/6 + 1.4 + 0.75i)$ of Example~\ref{exmp:hypatia-skew}, shown are 
box chain components in four fibers: left to right over $-10=-\beta, -9=\alpha, 9=-\alpha, 10=\beta$, with boxes from a $(2^{9})^2$ $w$-grid on $[-2.426, 2.426]^2$.  
}
\label{fig:matrab}
\end{figure}

\end{exmp}

%==============================================================
\subsection{Computability of $\mathcal{R}$ for polynomial skew products}
\label{sec:proof-thm-1}

Let $f:\bC^2\to\bC^2$ be an Axiom A polynomial skew product.
To establish the computability of $\mathcal{R}_f=\Omega_f=\overline{\text{Per}(f)}$ for $f$, by Corollary~\ref{cor:set-computable}, we construct an algorithm that on input $N\in\bN$ 
computes a set which is a $2^N$-approximation of $\cR$; i.e., within the $2^{-N}$-neighborhood of $\cR$ in the Hausdorff metric.

For any positive integer $n$, Algorithm~\ref{alg:box-cover} computes a set of boxes $\mathcal{B}_n$ such that
$\mathcal{R} \subset  \mathcal{B}_n \subset \mathcal{R}(\alpha),$ for any $\alpha > \ep'_n+2\xi_n$ (see Corollary~\ref{cor:pseudo-orb-calc}).

Since the boxes are defined in terms of pseudo-orbits, $\cB_n$ is obviously not a neighborhood of $\mathcal{R}$ in the Hausdorff metric, so it doesn't immediately imply computability. However, since the algorithm produces a nested sequence of sets whose intersection is $\mathcal{R}$ (Theorem~\ref{thm:boxchainconstruction}), it is reasonable to explore whether this nested sequence can be used to establish computability.

Recall that if the map $f$ is Axiom A, each chain component has what we call a ``hyperbolic type''; it is either expanding, $J_2$-type, saddle type, or attracting type. For saddle-type components, the attraction is either in the base with vertical expansion in the fibers, or the base is the Julia set of the one-dimensional polynomial, and the attraction is in the fibers.

To compute the chain components, we first show that if $f$ is Axiom A, for sufficiently large $n$ the box chain components of Algorithm~\ref{alg:box-cover} separate invariant subsets of different hyperbolic types, and we provide an algorithm to identify the hyperbolic type of each.
Note that if $n$ is not sufficiently large, there may be ``fake'' components; we shall deal with those in a subsequent algorithm.  Heuristically, one can quickly identify the box chain component that contains $J_2$ because it is the largest by far (i.e., has the most vertices and edges).

\begin{alg}\textbf{Computing a box chain model for a polynomial skew product such that every component is of hyperbolic type.}
\label{alg:compute-R-firststep}

(1) \textbf{Choose} a sequence of rationals $\gamma_k \downarrow 0$ as $k\to\infty$ (e.g.,$\gamma_k = 1/2^k$). 

(2) \textbf{Choose} a sequence of integers $\nu_k \uparrow \infty$,  (e.g., $\nu_k = 2^k$). 

(3)
For ease of exposition, set $g=f^{\nu_k}$. 

Note that $g$ is also a polynomial skew product, just of higher degree than $f$. 
For $f(z,w)=(p(z),q(z,w))$ and $x=(z,w)\in \bC^2$  write $z_i = p^i(z)$ then $Q^i_z = q_{z_{i-1}} \circ \cdots \circ q_z,$ so that $f^i(x) = f^i((z,w)) = (p^i(z),Q^i_z(w)). $
Next, recall that $Dg$ has a computable modulus of continuity
as described in the proof of Lemma~\ref{lem:distance-expansion_cw}. Thus, we can \textbf{compute} a power of $1/{2}$, denoted by $\Delta_k$,  which is small enough that two points in $\cB_0$ who are at most $\Delta_k$ apart have their derivatives differing by less than $\gamma_k$. In fact, we require that  $\inorm{D_xg}, \supnorm{D_xg}, |p^{\nu_k}(z)|,$ and $\onorm{DQ_z^{\nu_k}(w)}$ differ by less than $\gamma_k$ for points at most $\Delta_k$ apart. 
Thus, if $x=(z,w)$ and $x'=(z',w')$, we have 
\begin{align}\label{4ineq}
\norm{x-x'}\leq \Delta_k \Rightarrow & &   \\
& \inorm{D_xf^{\nu_k} - D_{x'}f^{\nu_k}}  &< \gamma_k,\\
& \supnorm{D_xf^{\nu_k} - D_{x'}f^{\nu_k}} &< \gamma_k, \\
& \onorm{p^{\nu_k}(z)-p^{\nu_k}(z')} &< \gamma_k, \\
& \onorm{DQ_z^{\nu_k}(w) -DQ_{z'}^{\nu_k}(w')} &< \gamma_k.
 \end{align}

Recall that the boxes in $\mathcal{B}_n$ computed by Algorithm~\ref{alg:box-cover} are of sidelength $\ep_n=R/2^{n-1}=2^{m-n+1}$, where $R=2^m$ depends on $f$. Thus, we can calculate $n=n_k=n(\nu_k,\gamma_k)$ so large that $ 2^{m-n+1} \leq \max(1, \Delta_k$).
Thus, for every box $B \in \mathcal{B}_n$, for any $x,x'\in B,$ the four inequalities above hold.

Now \textbf{calculate}  $n=n_k$ large enough and \textbf{run} the box chain construction of Algorithm~\ref{alg:box-cover}, using the function $g$.  

Finally,  \textbf{decompose} $\cB_n$ into its finite number of box chain components $\cB_n =  \cB_{n,1} \sqcup \ldots \sqcup \cB_{n,s}.$

(4) Perform a loop considering one at a time each box chain component $\cB_{n,i} \in \{\cB_{n,1}, \ldots, \cB_{n,s}\}$ found in Step (3). For a given $\cB_{n,i}$,
 for each box $B_l^n \in \mathcal{B}_{n,i}$, denote the center point of the box $B^n_l$ by $x^n_l=(z^n_l,w^n_l)$.

(a) [\textit{expansion: $J_2$-type}] 
For each  $B_l^n \in \mathcal{B}_{n,i}$,
{\bf test} if
$ \inorm{D_{x^n_l} f^{\nu_k}}  > 1 + \gamma_k$.
Note $\norm{D_x f^{\nu_k}}_{\inf}$ can be computed at any precision because it is the reciprocal of the sup norm of the matrix inverse.

If this statement fails for any box, stop and move on to step (b) below. 

Otherwise, this statement is true for every box in this component. {\bf Label} this component to be $J_2$-type, expanding. {\bf Set} $\lambda_i =  \min_{B_l^n \in \cB_{n,i}} \inorm{D_{x_l^n}f^{\nu_k}}-1-\gamma_k$, so that 
  $\inorm{D_{x}f^{\nu_k}} \geq 1+\lambda_i, \forall x\in\cB_{n,i}.$ Note $\lambda_i>0$ since there are a finite number of boxes.
Then, go back to the start of step (4), replacing $i$ with $i+1$ and considering the next component. 

(b) [\textit{contraction: $X_0$-type}] For each $B_l^n \in \mathcal{B}_{n,i}$, {\bf test} if
$ \supnorm{D_{x^n_l} f^{\nu_k}}  < 1 - \gamma_k$.
If this statement fails for any box, stop at that box, then move on to step (c) below. 

Otherwise, this statement is true for every box in this component. Then {\bf label} this component to be of $X_0$-type, attracting, and {\bf set}
$\lambda_i = 1 - \gamma_k- \max_{B_l^n \in \cB_{n,i}} \supnorm{D_{x_l^n}f^{\nu_k}}.$
So, $\supnorm{D_{x}f^{\nu_k}}  \leq 1-\lambda_i$, $\forall x  \in \cB_{n,i}$, with again $\lambda_i>0$. 
Then, go back to the start of step (4), replacing $i$ with $i+1$ and considering the next component. 

(c) [\textit{$z$-contraction, $w$-expansion: $X_1$-type}]
For each $B_l^n \in \mathcal{B}_{n,i}$ {\bf test} if both
$|p^{\nu_k}(z^n_l)|  < 1 - \gamma_k$
and 
$ \onorm{DQ_{z^n_l}^{\nu_k}(w_l^n)}  > 1 + \gamma_k$ are true.

If either statement fails for any box, stop at with that box, and move on to step (d) below.

Otherwise, both statements are true for every box in this component. {\bf Label} this component to be the $X_1$-type, saddle over $A_p$, and 
set 
$\lambda_i = \min \{ 1-\gamma_k-\max_{B_l^n\in \cB_{n,i}} | p(z_l^n)| , \min_{B\in \cB_{n,i}}|DQ^{\nu_k}_{z_l^n}(w_l^n)|-1-\gamma_k\}.$
So, we have $|p^{\nu_k}(z)|  < 1 - \lambda_i$ and $ \onorm{DQ_{z}^{\nu_k}(w)}  > 1+\lambda_i$ for all $(z,w)\in \cB_{n,i}.$
Then, go back to the start of step (4), replacing $i$ with $i+1$ and considering the next component. 

(d) [\textit{$z$-expansion, $w$-contraction: $\Lambda$-type}]
For each $B_l^n \in \mathcal{B}_{n,i}$, {\bf test} whether both
$|p^{\nu_k}(z^n_l)|  > 1 + \gamma_k$
and 
$ \onorm{DQ_{z^n_l}^{\nu_k}(w_l^n)}  < 1 - \gamma_k$ hold.

If either statement fails for any box, stop at that box and move on to step (e) below.

Otherwise, both statements hold for every box in this component. Then {\bf label} this component to be of the $\Lambda$-type, saddle over $J_p$, and
{\bf set} $\lambda_i = \min \{ \min_{B_l^n\in \cB_{n,i}}|p(z_l^n)|-1-\gamma_k, \ 1-\gamma_k-\max_{B_l^n\in \cB_{n,i}} |DQ^{\nu_k}_{z_l^n} (w_l^n) | \},$
so we have $|p^{\nu_k}(z)|  > 1 + \lambda_i$ and $ \onorm{DQ_{z}^{\nu_k}(w)}  < 1-\lambda_i$ for all $(z,w)\in \cB_{n,i}.$
Then, go back to the start of step (4), replacing $i$ with $i+1$ and considering the next component.

\smallskip

(e) If we reach this step, then we have a component without a hyperbolic type. In that case, stop this Step (4) loop (do not continue to look at the components for this $n$), then start over with Step (3) but replacing $k$ with $k+1$, and thus increasing $n=n_k$, resulting in smaller boxes. 
\qed
\end{alg}

If for some $k$, we don't hit step (4e), then the algorithm terminates and has produced a collection $\cB_n = \cB_{n,1}\sqcup  \ldots\sqcup  \cB_{n,s}$, where for each component $\cB_{n,i}$ we know the hyperbolic type ($J_2, X_0, X_1,$ or $\Lambda$) and we have a  $\lambda_i$ as defined above. We call this $\lambda_i$ the ``cushion of hyperbolicity'' for $\cB_{n,i}$. 

Next, we show that the algorithm above terminates if $f$ is Axiom A. 

\begin{lem} \label{lem:halting}
Suppose $f$ is an Axiom A polynomial skew product of $\CTwo$. Then the main loop of Algorithm~\ref{alg:compute-R-firststep} terminates for some $n\in\bN$.     
\end{lem}

\begin{proof}
We claim that for sufficiently large $n$, i.e., sufficiently small box size, the main loop of Algorithm~\ref{alg:compute-R-firststep} determines each component of $\cB_n$ to be one of the hyperbolic types, and thus halts. This follows from the following observations:

First, the chain recurrent set $\cR$ is the intersection of the decreasing sequence of the compact sets $\overline{\cR_{\alpha}}$, because (i) $\cR_\alpha \subset \cR_{\alpha'}$ if $\alpha \leq \alpha'$, (ii) $\cR = \cap_{\alpha>0} \cR_\alpha$, and (iii) each $\overline{\cR_{\alpha}}$ is compact. Hence, $d_H(\cR, \cR_\alpha)\to 0$ as $\alpha\to 0$. Recall that each component of a hyperbolic type is a finite union of basic sets all of which have a positive distance from each other. 
Since the nested sequence $(\cB_n)_{n\in\bN}$ decreases to $\cR$ (see Theorem~\ref{thm:boxchainconstruction}, (4)), it follows that for $n$ large enough, the boxes are small enough so that the box chain components separate invariant subsets of different hyperbolic type. 

For a polynomial skew product $f$, we know that at least $ J_2 \neq \emptyset$. 
If $f$ is Axiom A, then in a neighborhood $\mathcal{N}$ of $J_2$ there is a metric, uniformly equivalent to the Euclidean metric, for which $f$ is expanding; hence it is eventually expanding with respect to the Euclidean metric. That is, there is an iterate $\nu\in\bN$ and $\lambda > 0$ such that for all $x \in \mathcal{N}$ and $j \geq \nu$, 
we have $\inorm{D_xf^{j}} > 1+\lambda$.  
Since Algorithm~\ref{alg:box-cover} produces decreasing approximations of $J_2$  (the sequence $\cB_n$ tends to $\mathcal{R}$ as $n\to\infty$), for large enough $n$, some box chain transitive component $\cB'_n$ satisfies $J_2 \subset \mathcal{B}'_n \subset \mathcal{N}$. Since $\nu_k \uparrow\infty$, for $k$ large enough, $\nu_k>\nu$.
Since $\gamma_k \downarrow 0$, for $k$ large enough, $\lambda - \gamma_k > 0$. Therefore, for large enough $k$, both these statements hold. It follows that $\inorm{D_xf^{j}} > 1 + \gamma_k$ holds for every point in the box, and  hence in particular for the center points of all boxes. Thus, for large $n$, the algorithm above determines one of the box chain components to be of $J_2$-type. 

Since $f$ is an Axiom A polynomial skew product, any splitting of the tangent bundle of a basic set is horizontal/vertical (i.e., in $z$ respectively $w$ direction), and there are finitely many basic sets of $f$. Analogously to the $J_2$ argument, each basic set is of one of the remaining three hyperbolic types and thus each has a neighborhood around it in which the norms of the derivatives display the hyperbolic type, as described in the $J_2$ case above. Since there are finitely many basic sets, there is a common lower bound on the neighborhood size say $\delta$, so that in $\cN(\cR,\delta)$, the $\delta$-neighborhood of $\cR$,  all points show the hyperbolic type of the basic sets.
Since  $\cB_n\downarrow \cR$ as $n\to\infty$, there exists $n\in\bN$ (sufficiently large) that all of the box chain transitive components $\cB_n$ lie within $\cN(\cR,\delta)$. Hence each box chain transitive component of $\cB_n$ has a hyperbolic type detectable by the above algorithm, and we conclude that the algorithm halts.  
\end{proof}

We note that in Lemma \ref{lem:halting} above, it is possible that the algorithm produces box chain components which do not contain a basic set of a specific hyperbolic type; the invariant set $J_2$ is the only one which is guaranteed to be nonempty. This issue will be handled by the algorithm below, which is the final piece of our process for computing each type of hyperbolic invariant set. For each box chain component produced by Algorithm~\ref{alg:compute-R-firststep}, we perform an algorithm, tailored to the hyperbolic type of the component,
to refine the box chain component until it is an approximation of the hyperbolic invariant subset (in case it exists) to the desired precision.

\begin{alg}\textbf{for a given $N$, computing a $2^N$-approximation of the sets $J_2$, $X_0, X_1$ and $\Lambda$ of an Axiom A polynomial skew product.}

\label{alg:compute-R-increaseprecision}

Given a polynomial skew product $f$ of degree $d\geq 2$, perform Algorithm~\ref{alg:compute-R-firststep} to 
produce a collection $\cB_n = \cB_{n,1}\sqcup  \ldots\sqcup  \cB_{n,s}$, where for each component $\cB_{n,i}$ we know the hyperbolic type  ($J_2$, $X_0$, $X_1$ or $\Lambda$) and we have a  $\lambda_i$ as defined above for some iterate $g=f^{\nu_i}$, the cushion of hyperbolicity for $\cB_{n,i}$. 

Given $N\in\bN$ such that the $2^N$-approximation of $\cR(f)$ has these desired properties, perform one of the following four algorithms for each box chain component, based on the hyperbolic type of the component. 

Or, if a $2^N$-approximation is not desired for all of $\cR(f)$, but simply for $J_2,X_0, X_1$ and/or $\Lambda$, perform the appropriate algorithm below for the component(s) of the desired type. 

For ease of notation, we write $\cB'_n$ instead of $\cB_{n,i}$ for the box chain component under consideration in each individual algorithm below. Since each algorithm is a separate subroutine, we may re-use ``local variables''; that is, a constant introduced in one algorithm may be re-used in a different way in another one. 

\smallskip

\underline{(a) Type $J_2$ algorithm}:
Actually, since we know there is exactly one invariant set $J_2$ of the $J_2$-type, rather than considering one $J_2$-type box chain component below, let $\cB'_n$ in this sub-algorithm denote the union of all $J_2$-type box chain components in $\cB_n$. (There could be more than one, as there could be a ``fake'' expanding component near $J_2$ for some $n$.)

We start by outlining our approach. 
 We begin with the collection of all boxes $\cB'_n$ produced by Algorithm~\ref{alg:compute-R-firststep}, of $J_2$-type, and on which we have expansion by $Dg$ with cushion of hyperbolicity $\lambda'$. 
We apply Algorithm~\ref{alg:box-cover} to refine this set to a subset of smaller boxes $\cB'_{n+t}$ (which is similarly the union of all box chain transitive components which are subsets of $\cB'_n$), 
where the size of these smaller boxes is calculated to be so small that all points in $\cB'_{n+t}$ are $\alpha$-pseudo-periodic points (with $\alpha$ related to the box size in $\cB'_{n+t}$ and bounds on the derivative of the map), where we use a result from Urbanski et al.\ \cite{Urbanski_book}'s shadowing calculations to find the required upper bound on $\alpha$ needed, based on $\lambda'$ and $\beta_N:=2^{-N}$, to guarantee that all $\alpha$-pseudo periodic points in this expanding set are $\beta_N$-shadowed by a true orbit in $J_2$. Thus each point in $\cB'_{n+t}$ is within $\beta_N=2^{-N}$ of a point in $J_2$. 

We now present the proof in more detail.
First, suppose $\cB'_n$ is all box chain components produced by Algorithm~\ref{alg:compute-R-firststep} and determined to be of $J_2$-Type, uniformly expanding. 
Thus, for all points in $\cB_n'$,  we have $\inorm{D_{(z,w)}g}\geq 1+\lambda'$, where $g=f^{\nu_k}$, and $\lambda'$ is the cushion of hyperbolicity for $\cB'_n.$

(1) Using Corollary~\ref{cor:distance-expansion-cw} applied to $g$ on $\cB'_n$, we compute the radius $r'$ so that if two points in $\cB'_n$ are in a box of radius $r'$, then their images under $g$ are pushed apart by at least $1+\lambda'$. By taking $r'$ slightly smaller, assume it's a dyadic rational. 

(2) Next, first increase $N$ if needed (to obtain a more accurate approximation that's acceptable), so that $2^{-N}\leq 2r'$.

(3) Recall from Theorem~\ref{thm:boxchainconstruction}, the sequence $\{\xi_i\}$ is chosen by the user and converges to $0$. From Corollary~\ref{cor:imageboxsize}, $\ep'_i = (1+L_i)\ep_i$ where $L_i$ is computable and nonincreasing, and $\ep_i = 2^{m-i+1},$ given the original bounding box is radius $R=2^m$. 
Hence, we can compute a  positive integer $t_N$ sufficiently large that for all $t\geq t_N$: 
$$
\ep'_{n+t}+2\xi_{n+t} = (1+L_{n+t})\ep_{n+t} + 2\xi_{n+t} <  \min \{2r', \lambda' (2^{N}) /2 \}.
$$
Note this inequality guarantees $t\geq t_N$ is so large that $\ep_{n+t} \leq 2r'$. We see below why this inequality is desired.  

(4) Now for any $t\geq t_N$, apply Algorithm~\ref{alg:box-cover} to refine $\cB'_n$ by dividing each box into a $(2^t)^4$-grid.  That is, just apply the algorithm with $\cB'_n$ as the initial set, and let $\cB'_{n+t}$ be the union of all box chain components whose boxes are subsets of boxes in $\cB'_n$ (though we know only one of them contains $J_2$).
The boxes in $\cB'_{n+t}$ are of sidelength $\ep_{n+t}=2^{m-n-t+1}\leq 2r'$. 

Since the $\cB_n$ sequence is nested, points in some box in $\cB_{n+t}$ are also in some (larger) box $\cB_n$. So, we still have (at least) the same expansion, and by choice of the max box sidelength $\leq 2r'$, we have that $g$ is distance expanding by $1+\lambda'$ in each box in $\cB'_{n+t}$. 

(5) We claim 
$d_H(J_2, \cB'_{n+t}) \leq 2^{-N}$. 

We use shadowing to find how close a real orbit is to our pseudo-orbits. 
Since $g$ is  Axiom A (because $f$ is), it has the shadowing property on each chain component.
We apply Equation (4.29) in \cite{Urbanski_book} to our setting. This provides a minimum (denoted by $\xi$ in \cite{Urbanski_book}, but our $\xi_n$ is not the same) which in our setting is $2r'$ because the largest ball that fits inside of an image of a ball of radius $2r'$ has a radius bounded by $\inorm{Dg} 2r' 
\geq (1+\lambda')2r' > 2r'$ since $\inorm{Dg}\geq 1+\lambda' > 1$. 
By Proposition 4.3.4 in \cite{Urbanski_book}, we deduce that $\beta < 2r'$, and for $\alpha = \min (2r', \lambda'\beta/2))$, there is a (unique) true orbit in $J_2$ that $\beta$-shadows any $\alpha$-pseudo-orbit. 
(We note that if the expansion satisfies $\lambda' \geq 4r'/\beta,$ we have $\alpha=2r'$. On the other hand, if expansion satisfies $\lambda < 4r'/\beta,$ we obtain $\alpha < 2r'$.)

To translate this back to box sizes, in any $\cB'_{n+t}$ with boxes of size $\ep_{n+t}=2^{m-n-t+1}\leq 2r',$  we know by Corollary~\ref{cor:pseudo-orb-calc} that the points in boxes of size $\ep_{n+t}$, given a $\xi_{n+t}$ error bound in calculating the images of boxes under $f$, are in pseudo-orbits of size $> \ep'_{n+t}+2\xi_{n+t}$.

Recall $2^{-N}\leq 2r'$ and consider $\beta_{N} :=2^{-N}$.
Now $\beta_N$ defines $\alpha = \alpha(\beta_N) = \min (2r', \lambda'\beta_N/2 )$, for which $\alpha(\beta_N)$-pseudo-orbits are $\beta_N$-shadowed by true orbits, and these are orbits in $J_2$ by choice of the chain transitive component. 

By choice of $t_N$, we conclude that
the box size $\ep_{n+t}=2^{m-n-t+1}$ in $\cB'_{n+t}$ is sufficiently small that 
$ \ep'_{n+t}+2\xi_{n+t} < \min (2r', \lambda'\beta_N/2 )=\alpha(\beta_N).$

Thus, we computed a $t_N$ sufficiently large that for all $t\geq t_N$, any point in any box in $\cB'_{n+t}$ is in an $ \alpha(\beta_N)$-pseudo orbit, and thus is (strictly) within $\beta_N=2^{-N}$ of a true orbit. Hence $\cB'_{n+t}$ lies within a $\beta_N=2^{-N}$ open neighborhood of $J_2$, i.e., $\cB'_{n+t} \subset \cN(J_2, 2^{-N})$.

But we also know $J_2 \subset \cB'_{n+t}$.
So, we have 
$$
J_2 \subset \cB'_{n+t} \subset \cN(J_2,2^{-N}),
$$ 
for any $t \geq t_N$.
Hence $d_H(J_2,\cB'_{n+t}) \leq 2^{-N},$ if $t \geq t_N$. 
Thus $\cB'_{n+t}$ is a $2^{N}$-approximation of $J_2$.

\medskip
 This establishes the computability of $J_2$. We use variations of this algorithm to establish the computability of some of the remaining hyperbolic components.

We start by considering the two types of hyperbolic components (attracting and saddle-type) where the hyperbolic components lie in fibers over (finitely many) attracting periodic points $A_p$ in the base. In this situation, the computability problem becomes one-dimensional. To see this, let $f(z,w)=(p(z),q(z,w))$ be an Axiom A skew product and let $z$ be an attracting periodic point with period $\eta$ of the base polynomial $p$. Recall that since $p$ is hyperbolic there are only finitely many such points $z$. Considering $f^\eta$ instead of $f$ we may assume that $z$ is an attracting fixed point. In this case (for $z$ being fixed), the  $w$-component of $f^\eta$ is the one-dimensional polynomial $Q^\eta_z(w) =  q_{z_{\eta-1}}\circ \cdots \circ q_z$. The problem now reduces to computing the attracting cycles of $Q^\eta_z$ (attracting component) and the Julia set of $Q^\eta_z$ (saddle component).\\[0.01cm]

\underline{(b) Type $X_0$ algorithm}: Now, suppose $\cB'_n$ is a box chain component produced by Algorithm~\ref{alg:compute-R-firststep} and determined to be of Type $X_0$, an attractor. 

In this case, the base attracting set is a finite union of attraction cycles (periodic orbits) which can be computed at any given precision by one of the standard algorithms using that the basin of attraction of an attracting cycle contains a critical point. Here, we also use the standard fact that the critical points of a polynomial can be computed at any precision. We now consider for each periodic point the map $f^\eta$ instead of $f$, where $\eta$ is a period of the attracting period point. For this iterate, the fiber map becomes a one-dimensional polynomial. Thus, we can compute the attracting cycles of the fiber polynomial with the same algorithm as in the base case.
We note that these computations with sufficient accuracy may reveal this particular $\cB'_n$ as a ``fake'' component, which helps us to find all box chain components which lie within a given precision neighborhood of any global attractors. 

\smallskip

\underline{(c) Type $X_1$ algorithm}: Now suppose $\cB'_n$ is a box chain component produced by Algorithm~\ref{alg:compute-R-firststep} and determined to be of Type $X_1$: $z$-contraction, $w$-expansion.

We start with the same algorithm as in the Type $X_0$ case to compute the finitely many base-attracting periodic points at any given precision. 
This saddle set must lie over an attracting cycle in the base map. 
Replace the map with an iterate (as in Type $X_0$), to ensure that the fiber is fixed over the cycle, and within these fibers, the fiber Julia set is expanding. 
The Type $J_2$ algorithm can be simplified to one dimension and applied on that fiber. Note there are other algorithms to compute the fiber Julia set, e.g., \cite{Braverman2005}.

We also note that while this computation (with high enough accuracy) may show that this particular $\cB'_n$ is a ``fake'' component, it identifies all box chain components which lie within a given precision neighborhood of any Type $X_1$ saddle sets. 

\smallskip

\underline{(d) Type $\Lambda$ algorithm}:
Finally, suppose $\cB'_n$ is a box chain component produced by Algorithm~\ref{alg:compute-R-firststep} and determined to be of Type $\Lambda$: $z$-expansion, $w$-contraction. Thus, there exist $\lambda'>0$ and $\nu_k\in \bN$ such that for $g(z,w)=f^{\nu_k}(z,w)=(p^{\nu_k}(z),Q_{z}^{\nu_k}(w))$ we have $|p^{\nu_k}(z)|  > 1 + \lambda'$ and $ \onorm{DQ_{z}^{\nu_k}(w)}  < 1-\lambda'$ for all $(z,w)\in \cB_{n}'.$ Recall that any invariant set of this type lies in $J_p \times \CC.$

For clarity, we first outline the approach. 
We calculate an $\ep>0$ such that if we refine the boxes in $\cB'_n$ further, to diameter at most $\ep$, we have $\frac52\ep < 2^{-N}$, where we are trying to find the $2^{-N}$ approximation of a saddle invariant subset $\Lambda' \subset \cB'_n$ (if it exists). We consider boxes $B$ of this $\ep$-diameter and their ``doubled boxes'', we call $2B$, with the same center as $B$ but with double diameter $2\ep$ (see Figure~\ref{fig:doubledbox}). 
\begin{figure}
\includegraphics[width=0.65\textwidth]{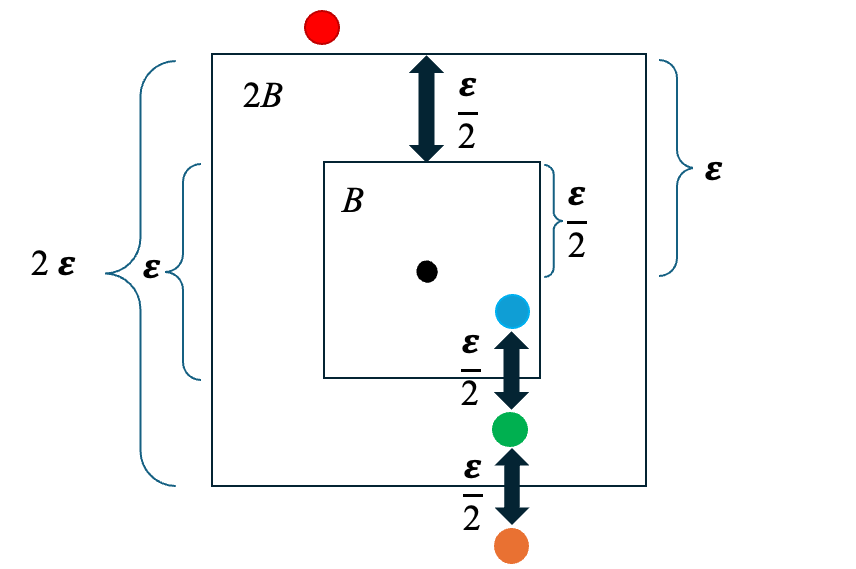}
\caption{\label{fig:doubledbox}
In the $\Lambda$-type algorithm, of Algorithm~\ref{alg:compute-R-increaseprecision}, we consider boxes $B$ of diameter $\ep$ and their ``doubled'' boxes with the same center but double diameter $2\ep$. We calculate approximate saddle periodic points of the map to accuracy $\ep/2$, and consider whether or not the saddle periodic points lie within $2B$.
}
\end{figure}
We run through a loop, indexed by an integer $j$ starting from $j=0$, each step refining the boxes further to a sidelength $\ep_j = \ep/2^j$. For each $j$, we first compute approximate repelling periodic points of the base map up to a sufficiently large period denoted by $u_j$ with $u_j > u_{j-1},$ such that every base box of diameter $\ep_j$ contains an approximate repelling periodic point of the base map. Then we compute approximations of \textit{all} saddle periodic points over each base repelling periodic point. We specify that both the repelling and the saddle periodic point calculations be accurate to within $\ep_j/2$, which is equal to the Hausdorff distance between a box $B$ of diameter $\ep_j$ and its double $2B$ (again, see Figure~\ref{fig:doubledbox}). Now, we examine each box $B$ in the level $j$ box collection. There are three cases for each $B$ denoted by Cases I, IIa, IIb. 

Case (I): The ``doubled'' box $2B$ may contain one of our approximate saddle points (the green dot in Figure~\ref{fig:doubledbox}). Since the accuracy of computing saddle points is $\ep_j/2$, it follows that either the box $B$ contains a saddle point (e.g., the blue dot in Figure~\ref{fig:doubledbox}), or a saddle point lies within $\ep_j/2$ of $2B$ (the green and orange dots in Figure~\ref{fig:doubledbox}). 
Thus, every point in the doubled box lies within $\frac52\ep = (2\ep + \ep/2)$ of a saddle periodic point. Since $\ep$ satisfies $\frac52\ep < 2^{-N}$, if this case occurs for \textit{all} boxes at step $j$, then the collection of all doubled boxes $2B$ forms a $2^{-N}$-approximation of the saddle invariant set contained in $\cB'_n$ (assuming such a set exists). 

Case (II): The doubled box $2B$ may not contain any of our approximate saddle points (they may all be contained in the complement of $2B$, like the red dot in Figure~\ref{fig:doubledbox}). Again, given our computational accuracy, because the distance from the complement of $2B$ to $B$ is $\ep_j/2$ (see Figure~\ref{fig:doubledbox}), this means that the box $B$ definitely does not contain any saddle periodic point of period smaller or equal to $u_j$. Since the saddle-type invariant sets are the closure of a set of saddle periodic points, either this means:
\\Case (IIa): the box $B$ truly contains no saddle periodic points of $g$, or 
\\Case (IIb): We just haven't checked large enough periods. 
\\At this point, we are not able to determine which case occurs. If we detect at loop $j$ at least one box which does not contain any of our computed saddle points, we stop the $j$-loop and start over with $j+1$, and perform the algorithm for smaller boxes and higher period saddle points.  

We stop increasing $j$ when all doubled boxes at that refinement level are in Case (I); i.e., they contain an approximate saddle point (of period at most $u_j$). We show in the proof of Theorem~\ref{thm:main} below that the algorithm terminates if the set of boxes $\cB'_n$ contains a hyperbolic set of type $\Lambda$. This essentially follows from the facts that a saddle-type hyperbolic set of $f$ is the closure of its saddle points and that our collection of boxes chosen and refined as described in Algorithm~\ref{alg:box-cover} is a nested sequence of compact sets that converges to the chain recurrent set. We note that it is in principle possible that this component is ``fake", and at some point in the $j$-th loop, a refinement contains no boxes in cycles in which case the algorithm terminates. This will be discussed in more detail below.

We now describe the detailed algorithm. 

(1) 
Recall the boxes in $\cB'_n$ are of diameter $\ep_n$ which is a power of $1/2$; in fact, $\ep_n = R/2^{n-1}=2^{m+1-n}$. 
Compute  $t\in\bN$  large enough that 
$\frac52\ep_{n+t} = \frac52 (2^{m+2-n-t}) \leq 2^{-N}$
holds.

Set $j=0$.

(2)  Subdivide the boxes in $\cB'_{n}$ by placing a $(2^{t+j})^4$-grid on each box, and perform Algorithm~\ref{alg:box-cover} on $\cB'_n$. This might result in multiple box chain components at level $n+t+j$ all inside of $\cB'_n$. 
Since it is not clear which box chain components (if any) contain at least one basic set, the algorithm is performed on all of them. Recall that $\cB_{n+t+j}$ is the collection of boxes that contain all hyperbolic invariant sets of all hyperbolic types, at level $n+t+j$. Let $\cB'_{n+t+j}$ be the union of all box chain components of $\cB_{n+t+j}$ which lie in $\cB'_n$, if any! There may be no boxes in $\cB_{n+t+j}\cap \cB'_n$. In this case, we conclude that $\cB'_n$ was a ``fake'' component and we stop the main algorithm, and go back to the beginning of Algorithm~\ref{alg:compute-R-increaseprecision} and move on to the next (if any) box chain component of $\cB_n$ which has not been computationally analyzed. 

Now, if we have $\cB'_{n+t+j}\neq \emptyset$, continue the algorithm as follows. 

(3) Let $2\cB'_{n+t+j}$ denote the collection of boxes with the same center points as the boxes in $\cB'_{n+t+j}$, but with double the diameters. 
For a box $B$ in $\cB'_{n+t+j}$, we use $2B$ to refer to the box of the form $\cN(B,\ep_{n+t+j}/2):$, i.e., $\ep_{n+t+j}/2$-neighborhood in our $L^\infty$ metric of the box $B$. (See Figure~\ref{fig:doubledbox}.)

(4) Since $B$ is a box in $\CTwo$ we can write $B=B_z \times B_w$, where $B_z$ and $B_w$ are one-complex dimensional boxes in the $z$-plane and  $w$-plane, respectively.  Let $(\cB_{n+t+j}')_z$ denote the collection of $z$-boxes $B_z$ of all boxes $B$ in $\cB'_{n+t+j}$.

Recursively compute a strictly increasing sequence $\{u_j\}_{j\in\bN}$ as follows.

Compute an approximation to all repelling periodic points of the base map $p^{\nu_k}(z)$ of $g$  (one period at a time) from period $1$ up to period $u_j$, where $u_j$ satisfies that every $B_z$  in $(\cB_{n+t+j}')_z$ contains contains at least one repelling periodic point of period larger than $u_{j-1}$ and smaller than $u_j$. 
This is possible because the Julia set of the base map is the closure of repelling periodic points of the base map.
Carry out the computation of these repelling periodic points with a precision of at least $\ep_{n+t+j}/2.$

(5) In this step, we perform a loop through each base box $B_z$ in $(\cB_{n+t+j}')_z$.

(i) For a base $B_z$ under consideration, we have identified repelling periodic points of $p^{\nu_k}$ contained in $B_z$. Let $\zeta$ be such a repelling period point and let $\eta$ be the period of $\zeta$. Then the $\eta$-th iterate of the fiber map of $g$, $(Q^{\nu_k})^\eta_\zeta$, is a one-dimensional complex polynomial depending on $\zeta$. 

Next compute approximations of \textit{all} attracting periodic points of  $(Q^{\nu_k})^\eta_\zeta$.  By increasing the accuracy of the computation of $\zeta$ if necessary, we can assure that the accuracy of the computation of these attracting periodic points is at least  $\ep_{n+t+j}/2$. This yields (finitely many) saddle points of $g$, and the accuracy of the computation of these saddle points is at least $ \ep_{n+t+j}/2$  in the $L^\infty$-metric.  

(ii) For the box $B_z$ there finitely many boxes in $\cB'_{n+t+j}$ that have $B_z$ as $z$-component. Next we perform a loop for each of these $\CTwo$ boxes $B$ in $B_z \times \CC$ to determine whether or not $2B$ 
contains one of the approximate saddle periodic points associated with $\zeta$.

If, while considering all boxes in $\cB'_{n+t+j}$ that have the $B_z$ under consideration as $z$-component, we detect that there is a $\CTwo$ box $B$ in $\cB'_{n+t+j}$ with the property that its doubled box does NOT contain any of computed saddle periodic points, then we stop this loop, 
and go back up to Step (2), proceeding with increasing $j$ to $j+1$.

We note that if the doubled box with sidelength $2\ep_{n+t+j}$  does not contain a saddle point that we calculated with accuracy $ \ep_{n+t+j}/2,$ it follows that the subset box of sidelength $\ep_{n+t+j}$ definitely does not contain a saddle periodic point with fiber period at most $u_j$, see Figure~\ref{fig:doubledbox}. 

Either a box $B$ contains no saddle periodic points, or it contains one of a higher period. This is why we increase $j$ and refine the box into a grid of smaller boxes and increase the base max period $u_j$.

On the other hand, if all doubled boxes over $B_z$ contain at least one approximate saddle periodic point, repeat (5)-(i) and (ii) for the next base box in the list of all base boxes at this level. 

Finally, if all doubled boxes for this level contain a computed saddle periodic point, move to procedure (6).

(6) If the algorithm reaches this step, we have examined all boxes in $\cB'_{n+t+j}$ and concluded that 
any point $x$ in $2\cB'_{n+t+j}$ lies in a box $2B$ of sidelength $2\ep_{n+t+j}$ which contains an approximate saddle point with accuracy $\ep_{n+t+j}/2$, hence there is true saddle periodic point within a distance of $\ep_{n+t+j}/2$ to the box $2B$. Hence $x$ is within the distance $(2\ep_{n+t+j}) + (\ep_{n+t+j}/2)=\frac52\ep_{n+t+j}$  from $\Lambda'$ (the closure of the set of saddle periodic points located in $\cB'_n$). Moreover, this distance  $\frac52\ep_{n+t+j}$ is smaller than $ 2^{-N}$ (by Step (1)). 
Thus, $2\cB'_{n+t+j} \subset \cN(\Lambda',2^{-N}).$ Since $\Lambda' \subset \cB_{n+t+j} \subset 2\cB'_{n+t+j}$, we have
$
\Lambda' \subset  2\cB'_{n+t+j} \subset \cN(\Lambda',2^{-N}),
$
thus $2\cB'_{n+t+j}$ is a $2^{-N}$-approximation of $\Lambda'$.
\qed

\medskip

Completing the appropriate type algorithm above for every box-chain component calculates $J_2, X_0, X_1$, and $\Lambda$ (if the set exists) to any desired precision. 
\qed

\end{alg}

\begin{proof}[Proof of Theorem~\ref{thm:main}]
With the algorithm above, we have proved the statement for the $J_2, X_0$, and $X_1$ type invariant sets. We only need to show that the algorithm for the $\Lambda$-type sets halts if $f$ is Axiom A.
This algorithm halts because the saddle periodic points are dense in $\Lambda$. If a box contains no saddle periodic point, we decrease the box size and refine. As we loop through smaller base boxes, we calculate saddle periodic points up to higher periods $u_j\to\infty$. If a non-empty invariant subset $\Lambda' \subset \cB'_n$ exists, 
the sequence of unions of collections of boxes $\cB'_{n+t+j}$ as described in the algorithm above is a compact, nested set tending to $\Lambda'\subset \Lambda$ as $j\to\infty$, and saddle periodic points are dense in $\Lambda$.  It may occur  that  $\cB'_{n+t+j}\cap \Omega_f =\emptyset$ holds. However, in this case, for sufficiently large $j$, the algorithm will terminate, since all boxes will disappear since all enlarged boxes do not contain pseudo orbits that are associated with orbits in the non-wandering set.
\end{proof}

%====================================
\section{Lower Computability of the Axiom A locus for polynomial skew products}
\label{sec:results2}
%====================================

We now present the proof of  Corollary~\ref{cor:AxiomAsemi-decidable}, i.e., that Axiom A is semi-decidable for polynomial skew products. 

\begin{proof}[Proof of Corollary~\ref{cor:AxiomAsemi-decidable}]
    If the Algorithm~\ref{alg:compute-R-firststep} does not reach Step (4e), it establishes bounds on the expansion and/or contraction of the individuel components, and it also establishes hyperbolicity on each component which contains an actual non-empty invariant set. Hence the algorithm detects Axiom A for an Axiom A polynomial skew product $f_c$ based on having oracle access to $c$,  as described in Lemma~\ref{lem:halting}. 
    
    Note that if $f$ lies on the boundary of the hyperbolicity locus, we won't be able to say it's not hyperbolic. 
    Since it's not known for complex polynomials whether hyperbolicity is dense in the entire parameter space, it can't be known for polynomial skew products. So the best we can hope for is semi-decidability. 
Thus, our method shows that Axiom A is a semi-decidable problem on the closure of the set of Axiom A polynomial skew products. 
\end{proof}

Next, we provide an algorithm to lower-compute the Axiom A locus, which proves  Theorem~\ref{thm:LowerComputabilityLocus}. 

We note that by (\cite{Jonsson1999}, Cor. 8.15) the set of Axiom A skew products on $\CTwo$ of fixed degree $d\geq 2$ is an open subset of the parameter space.

For fixed degree $d\geq 2$, recall $\mathscr{A}_d \subset \CC^{\ell}$ is the set of all Axiom A polynomial skew products of the form $f(z,w) = (p(z),q(z,w))$ where $p$ and $q$ are as in Equation~(\ref{eqn:coefficients_f}).

\begin{alg} \textbf{For a given degree $d\geq 2$, compute a countable list of ideal balls whose union is the locus of Axiom A polynomial skew products of degree $d$.}
\label{alg:compute-Hyp-locus}

(1) Let $M_k\uparrow \infty$ be positive integers and set  $\sM_k := [-2^{M_k},2^{M_k}]^{2\ell} \subset \CC^{\ell}$.

(2) Choose sequences $\gamma_k\downarrow0$ and $\nu_k\uparrow\infty$.  

(3) Let $\cS_k$ be the set of ideal points (points with dyadic rational coordinates) in $\sM_k$ of denominator $2^{-M_k}$. 
Each ideal point is a polynomial skew product $f_c$. 
Let $g_c$ denote $f_c^{\nu_k}$. 

(4) The idea is that for each $k$, run one loop of Algorithm~\ref{alg:compute-R-firststep} on all grid points, to effectively parallelize the algorithm. 

(i) For each $k$, run the algorithm just for that $k$: for that set of ideal points in $c \in \cS_k$, follow Algorithm~\ref{alg:compute-R-firststep} Step (3) to calculate the box collection based on this $k$, and then Step (4) to attempt to determine a hyperbolic type for each box chain component. 

(ii) For any ideal point map $f_c$ such that Algorithm~\ref{alg:compute-R-firststep} doesn't hit step (4e), but instead produces for each box chain component at this level a hyperbolic type, we know this map is Axiom A. If $\lambda$ is its cushion of hyperbolicity on its box chain components at this level, using the modulus of continuity of the derivative, similar to previous arguments but this time allowing the parameters of the polynomial (though not its degree just the coefficients) to change,
we can calculate a dyadic rational radius $r_c$ such that all maps in the neighborhood of radius $r_c$ in the parameter space are also Axiom A, with cushion of hyperbolicity at least $\lambda/2$. Add the ball of radius $r_c$ about this Axiom A ideal point $c$ to a list of balls whose union will be $\sA_d$. And for efficiency's sake, mark these ideal points, which are the centers of these balls as finished, so that on the next refinement (since dyadic rationals of denominator $2^{k+1}$ include the dyadic rationals of denominator $2^k$) we don't have to re-check any maps which we already determined were Axiom A.

Simply ignore any ideal points whose maps for this $k$ hit Step (4e). 

Since there are a finite number of ideal points in $\cS_k$, we add for each $k$ a finite number of balls to our union. 

(5) After running Steps (3)--(4) of Algorithm~\ref{alg:compute-R-firststep} on all ideal points in $\cS_k$, replace $k$ by $k+1$ and repeat Step (4) above, examining maps corresponding to ideal points of denominator $2^{-M_{k+1}}$ in $\sM_{k+1} = [-2^{M_{k+1}},2^{M_{k+1}}]^{2\ell}$.
\qed
\end{alg}

\medskip

\begin{proof}[Proof of Theorem~\ref{thm:LowerComputabilityLocus}]

We must simply argue why the above algorithm produces a union equal to $\sA_d$. 

Suppose a polynomial skew product $f_c'$ of degree $d$ is Axiom A. We explain why it is in one of the balls we have constructed. If $f_{c'}$ is Axiom A, by Lemma~\ref{lem:halting} there is a $k$ for which Algorithm~\ref{alg:compute-R-firststep} would produce hyperbolic type box chain components, if we ran the Algorithm on $f_{c'}$. But also, there is a radius $r_{c'}$ about $f_{c'}$ about which all degree $d$ polynomial skew products are Axiom A with perhaps a slightly smaller cushion of hyperbolicity than $f_{c'}$ has; not only are they Axiom A but they are in the locus of structural stability, the same hyperbolic component in parameter space as $f_{c'}$. There is also a sequence of ideal points in parameter space tending to $f_{c'}$, each with its own ball of a certain radius added to the union we compute in Algorithm~\ref{alg:compute-Hyp-locus}. Thus, there is an ideal point (infinitely many in fact) within this same hyperbolic component as $f_{c'}$, and for a close enough ideal point to $f_{c'}$, its ball of Axiom A maps which we added to our union will contain $f_{c'}$, because the radii were all constructed in a uniform way, using the modulus of continuity of the derivative maps as the coefficients of the polynomial vary. 
\end{proof}

%------------------------------------------------------
\subsection*{Future work}

This paper initiates the study of computability for holomorphic maps of $\CTwo$.  
One way to improve upon the results of this paper would be by providing a polynomial-time algorithm for computing the chain recurrent set. The given algorithm could be a first step, and a new algorithm could be provided to increase the precision. All current polynomial-time algorithms in one dimension use some type of conformality argument, which is not available in higher dimensions. This suggests new techniques are needed to control the computational complexity for maps in higher dimensions.

Polynomial skew products are a good entry point to the study of maps of higher dimension, as the splitting of the tangent bundle into stable and unstable directions is trivial. Thus, in future studies, one might examine complex \Henon diffeomorphisms, where the splitting is more complicated, with the splitting of the tangent bundle varying continuously and, in particular, moving around.

\appendix
\section{Examples of Axiom A polynomial skew products.} 
\label{sec:App_Examples}

\begin{exmp} \label{exmp:product}
The simplest case of an Axiom A skew product of $\bC^2$ is a \textit{product}, namely $f(z,w) = (p(z),q(w))$ is Axiom A if $p$ and $q$ are hyperbolic. There are (at most) four chain transitive components: $\cR = (J_p \times J_q)  \cup (A_p \times J_q)  \cup \ (J_p \times A_q)  \cup \ (A_p \times A_q)$; the first is $J_2$, the codimension zero set on which $f$ is uniformly expanding, the middle two are codimension one saddle sets, and the last set is the attracting periodic points of $f$.  For example, for $(z,w) \mapsto (z^2, w^2)$, we have $J_p \times J_q = S^1 \times S^1$ is a torus, $A_p \times J_q = \{0\} \times S^1$ and $J_p \times A_q = S^1 \times \{0\}$ are circles, and $A_p \times A_q$ is the origin. 
Moreover, by Jonsson's structural stability results, any $f(z,w)$ which is a sufficiently small perturbation of an Axiom A product is also Axiom A, with basic sets (i.e., chain components) topologically corresponding to those of the product.
\end{exmp}

In \cite{Boyd-DeMarco-2008, Boyd-DeMarco-2008E}, the first author and Laura DeMarco construct further examples of Axiom A polynomial skew products, including the following.

\begin{figure}
\includegraphics[width=.6\textwidth]{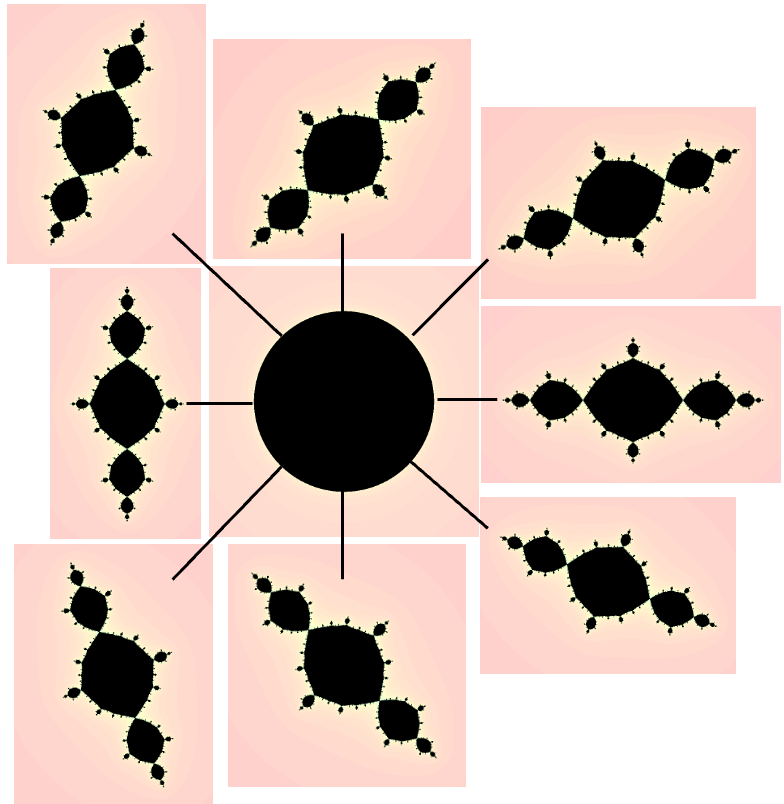}
\caption{\label{fig:example2}
$f(z,w) = (z^2, w^2 - z)$ with fiber $J_z$'s basilicas rotating as $z$ moves around the base $J_p$, which is the unit circle. Fibers $z=e^{i\theta}$ shown for $\theta$ values of $0, \pi/4, \pi/2, 3\pi/4, \pi, 5\pi/4, 3\pi/2,$ and $7\pi/4.$
$J_2$ is the boundary of the black region.}
\end{figure}

\begin{exmp} \label{exmp:twbas}
$F_a(z,w) = (z^2,w^2+az)$ is Axiom A iff $w \mapsto w^2+a$ is hyperbolic ($K_a$ can be connected or a Cantor set, either is permissible). This map is not in the same hyperbolic component as a product, though is semi-conjugate to a product via $(z,w) \mapsto (z^2, zw),$ at least on $\CTwo$. The fiber Julia sets are rotations of a quadratic Julia set. 
 See Figure~\ref{fig:example2}.
 \end{exmp}

\begin{figure}
\includegraphics[width=\textwidth]{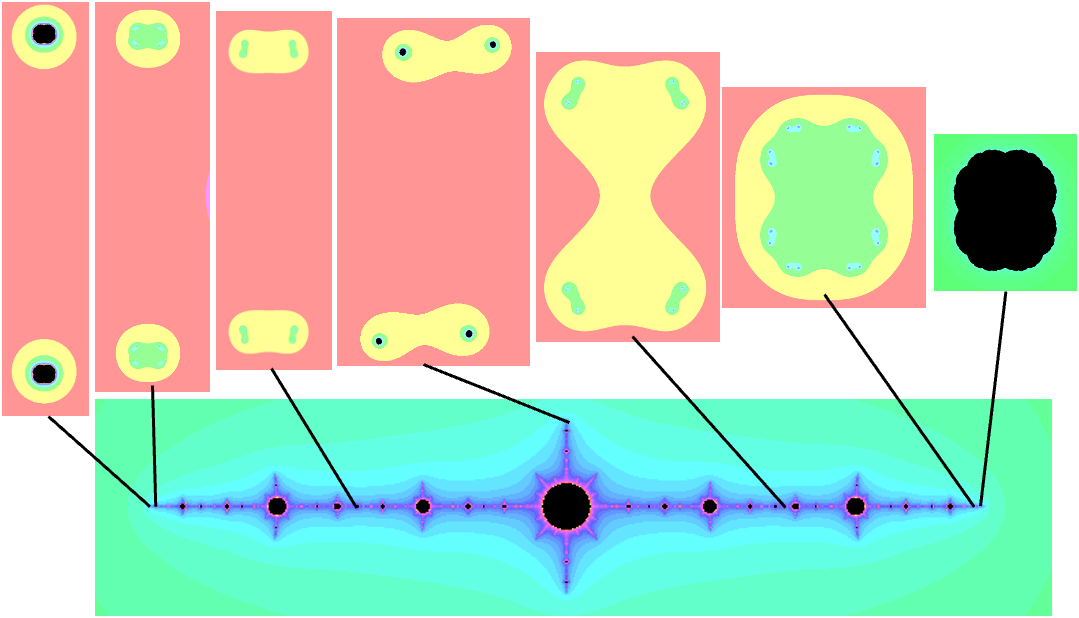}
\caption{\label{fig:example3}For $f(z,w) = (z^2-1.75488, w^2+2(2-z))$ we show $K_p$ (the aeroplane), 
and fibers for $z=-1.92, -1.8, -1, 0.4i, 1, 1.8, 1.92,$ 
with $J_z$'s a mixture of circles and Cantor sets.}
\end{figure}

\begin{exmp} \label{exmp:aerocantcirc}
Another type of example has an ``aeroplane'' base Julia set, $p_n(z) = z^2 + c_n$ is the unique quadratic polynomial with periodic critical point of lead period $n$ and with $c_n$ real.  The polynomial skew products 
$f_n(z,w) = (p_n(z), w^2 + 2(2-z))$ are Axiom A for sufficiently large $n$, and though $J_p$ is connected, we still have different fiber Julia sets over the two fixed points: one a quasi-circle and one a Cantor set.
See Figure~\ref{fig:example3}. 
\end{exmp}

\begin{figure}
\includegraphics[width=.95\textwidth]{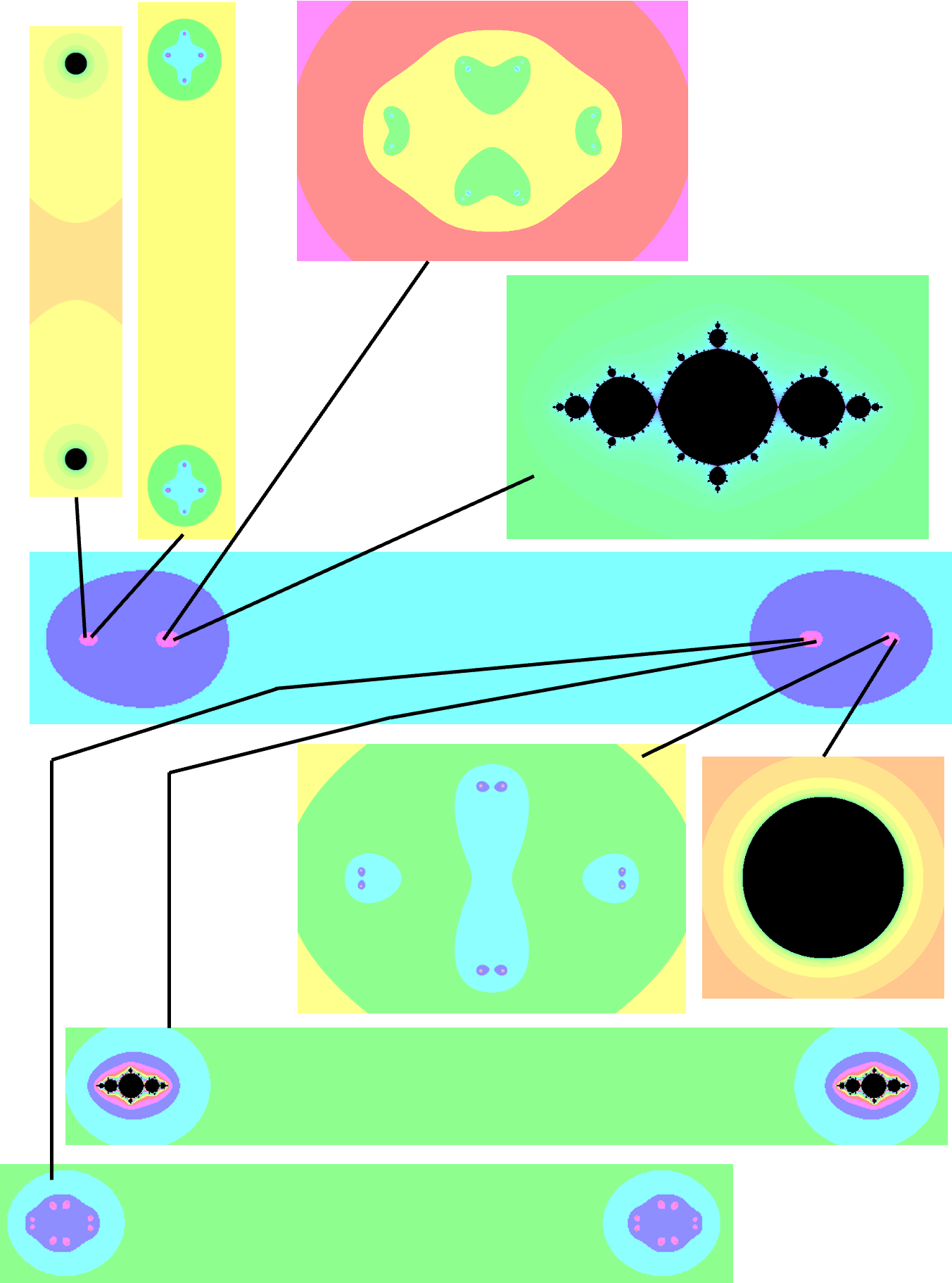}
\caption{\label{fig:example4} 
$f(z,w) = (z^2+20, w^2 + z^2 - 0.9z - 20.5)$, where one fixed fiber map is $w \mapsto w^2$, the other is $w \mapsto w^2-1$, i.e., circles and basillicas. 
We show $K_p$ (center of figure) and fibers for $z= -5, -4.99, -4.014, -4, 3.998, 4, 4.886,$ and $5$.  $K_5$ is the unit disk and $K_{-5}$ maps onto it, similarly $K_{-4}$ is a basilica and $K_4$ maps onto it. Other fibers are shown to have $J_z$ as a Cantor set.}
\end{figure}

\begin{exmp} \label{exmp:cantcircbas}
We provide a family of examples generalizing an example of Diller and Jonsson (\cite{MatDil}), in which the base is a Cantor set, and in the simplest case of degree two the fiber Julia sets over one fixed point are one hyperbolic Julia set, like a circle, and over the other fixed point a different one, say a basilica.
For example, see Figure~\ref{fig:example4}. 
\end{exmp}

Nekrashevych showed  the \textit{rational} skew product of $\Pt$ given by
$R(z,w) = \left(1 - 1/z^2, 1 - w^2/z^2 \right)$
is Axiom A, with a connected base Julia set and all fiber Julia sets connected, but not all fibers are homeomorphic. (For example, over the fixed points of the base map, one fiber map is the rabbit, while another one is the airplane).   This suggests there may exist polynomial skew products of $\CTwo$ with connected but varying fiber Julia sets, though we know of no such example.

%====================================

\clearpage
\bibliographystyle{plain}

\end{document}